\documentclass[twoside,a4paper]{article}

\usepackage{amsmath,amsthm,amssymb,latexsym}
\usepackage[v2,tips]{xy}

\newcommand{\noproof}{\hfill$\square$}
\newcommand{\wap}{\operatorname{WAP}}
\newcommand{\tobidual}[1]{\kappa_{#1}(#1)}

\newcommand{\ip}[2]{{\langle {#1} , {#2} \rangle}}
\newcommand{\aone}{\Box}
\newcommand{\atwo}{\Diamond}
\newcommand{\proten}{{\widehat{\otimes}}}
\newcommand{\mc}[1]{\mathcal{#1}}

\theoremstyle{plain}%
\newtheorem{proposition}{Proposition}[section]%
\newtheorem{theorem}[proposition]{Theorem}%
\newtheorem{corollary}[proposition]{Corollary}%
\newtheorem{lemma}[proposition]{Lemma}%
\theoremstyle{definition}%
\newtheorem{definition}[proposition]{Definition}%
\theoremstyle{remark}%
\newtheorem{example}[proposition]{Example}%

\begin{document}

\large
\title{\textsc{Connes-amenability of bidual and weighted semigroup algebras}}
\author{Matthew Daws}
\date{}
\maketitle

\begin{abstract}
We investigate the notion of Connes-amenability, introduced by
Runde in \cite{Runde1}, for bidual algebras and weighted semigroup
algebras.  We provide some simplifications to the notion of a
$\sigma WC$-virtual diagonal, as introduced in \cite{Runde2}, especially
in the case of the bidual of an Arens regular Banach algebra.  We apply
these results to discrete, weighted, weakly cancellative semigroup algebras,
showing that these behave in the same way as C$^*$-algebras with regards
Connes-amenability of the bidual algebra.  We also show that for
each one of these cancellative semigroup algebras $l^1(S,\omega)$, we have
that $l^1(S,\omega)$ is Connes-amenable (with respect to the canonical
predual $c_0(S)$) if and only if $l^1(S,\omega)$ is amenable, which is
in turn equivalent to $S$ being an amenable group.  This latter point
was first shown by Gr{\"o}nb\ae k in \cite{Gron1}, but we provide a
unified proof.  Finally, we consider the homological notion of
injectivity, and show that here, weighted semigroup algebras do
not behave like C$^*$-algebras.
\end{abstract}



\section{Introduction}

We first fix some notation, following \cite{Dales}.
For a Banach space $E$, we let $E'$ be its dual space, and
for $\mu\in E'$ and $x\in E$, we write $\ip{\mu}{x} = \mu(x)$
for notational convenience.  We then have the canonical
map $\kappa_E:E\rightarrow E''$ defined by
$\ip{\kappa_E(x)}{\mu} = \ip{\mu}{x}$ for $\mu\in E',x\in E$.
For Banach spaces $E$ and $F$,
we write $\mc B(E,F)$ for the Banach space of bounded linear
maps between $E$ and $F$.  We write $\mc B(E,E) = \mc B(E)$.
For $T\in\mc B(E,F)$, the \emph{adjoint} of $T$ is
$T'\in\mc B(F',E')$, defined by $\ip{T'(\mu)}{x} =
\ip{\mu}{T(x)}$, for $\mu\in F'$ and $x\in E$.

Let $\mc A$ be a Banach algebra.  A \emph{Banach left $\mc A$-module}
is a Banach space $E$ together with a bilinear map
$\mc A\times E \rightarrow E; (a,x) \mapsto a\cdot x$, such
that $\|a\cdot x\|\leq \|a\|\|x\|$ and $a\cdot(b\cdot x)=
ab\cdot x$ for $a,b\in\mc A$ and $x\in E$.  Similarly, we have the
notion of a \emph{Banach right $\mc A$-module} and a
\emph{Banach $\mc A$-bimodule}.  If $E$ is a Banach $\mc A$-bimodule
(resp. left or right module) then $\mc A'$ is a Banach
$\mc A$-bimodule (resp. right or left module) with module action
given by
\[ \ip{a\cdot\mu}{x} = \ip{\mu}{x\cdot a}
\qquad \ip{\mu\cdot a}{x} = \ip{\mu}{a\cdot x}
\qquad (a\in\mc A, x\in E). \]
Notice that as $\mc A$ is certainly a bimodule over itself (with
module action induced by the algebra product) we also have that
$\mc A'$, $\mc A''$ etc. are Banach $\mc A$-bimodules.  Given
a Banach $\mc A$-bimodule $E$, a subspace $F$ of $E$ is a
\emph{submodule} if $a\cdot x, x\cdot a\in F$ for each
$a\in\mc A$ and $x\in F$.  For Banach $\mc A$-bimodules
$E$ and $F$, $T\in\mc B(E,F)$ is an
\emph{$\mc A$-bimodule homomorphism} when
\[ a \cdot T(x) = T(a\cdot x) \qquad
T(x) \cdot a = T(x\cdot a) \qquad (a\in\mc A, x\in E). \]

A linear map $d:\mc A \rightarrow E$ between a Banach algebra $\mc A$
and a Banach $\mc A$-bimodule $E$ is a \emph{derivation} if
$d(ab) = a \cdot d(b) + d(a)\cdot b$ for $a,b\in\mc A$.  For $x\in E$,
we define $\delta_x : \mc A\rightarrow E$ by $\delta_x(a) =
a\cdot x-x\cdot a$.  Then $\delta_x$ is a derivation, called an
\emph{inner derivation}.

A Banach algebra $\mc A$ is said to be \emph{super-amenable} or
\emph{contractable} if every bounded derivation $d:\mc A\rightarrow E$,
for every Banach $\mc A$-bimodule $E$, is inner.  For example, a
C$^*$-algebra $\mc A$ is super-amenable if and only if $\mc A$ is
finite-dimensional.  It is conjectured that there are no infinite-dimensional,
super-amenable Banach algebras.

If we restrict to derivations to $E'$ for Banach $\mc A$-bimodules
$E$ then we arrive at the notion of \emph{amenability}.  For example,
a C$^*$-algebra $\mc A$ is amenable if and only if $\mc A$ is
nuclear; a group algebra $L^1(G)$ is amenable if and only if the
locally compact group $G$ is amenable (which is the motivating
example).  See \cite{RundeBook} for further discussions of
amenability and related notions.

Let $E$ be a Banach space and $F$ a closed subspace of $E$.
Then we naturally, isometrically, identify $F'$ with $E' / F^\circ$,
where
\[ F^\circ = \{ \mu\in E' : \ip{\mu}{x}=0 \ (x\in F) \}. \]

\begin{definition}\label{dual_means}
Let $E$ be a Banach space and $E_*$ be a closed subspace of
$E'$.  Let $\pi_{E_*}: E'' \rightarrow E''/E_*^\circ$ be the
quotient map, and suppose that $\pi_{E_*} \circ \kappa_E$
is an isomorphism from $E$ to $E_*'$.  Then we say that $E$
is a \emph{dual Banach space} with \emph{predual $E_*$}.

When $\mc A$ is a dual Banach space with predual $\mc A_*$ which
is also a submodule of $\mc A'$ we say that $\mc A$ is a
\emph{dual Banach algebra}.
\noproof
\end{definition}

For a dual Banach algebra $\mc A$ with predual $\mc A_*$, we
henceforth identify $\mc A$ with $\mc A_*'$.  Thus we get
a weak$^*$-topology on $\mc A$, which we denote by
$\sigma(\mc A,\mc A_*)$.  It is a simple exercise to show
that $\mc A$ is a dual Banach algebra if and only if
$\mc A$ is a dual Banach space such that the algebra product
is separately $\sigma(\mc A,\mc A_*)$-continuous (see
\cite{Runde1}).  The following lemma is standard.

\begin{lemma}\label{weak_star_cts}
Let $E$ and $F$ be dual Banach spaces with preduals $E_*$ and
$F_*$ respectively, and let $T\in\mc B(E,F)$.  Then the
following are equivalent:
\begin{enumerate}
\item $T$ is $\sigma(E,E_*) - \sigma(F,F_*)$ continuous;
\item $T'(\kappa_{F_*}(F_*)) \subseteq \kappa_{E_*}(E_*)$;
\item there exists $S\in\mc B(F_*,E_*)$ such that $S'=T$.
\end{enumerate}
\vspace{-2ex}\noproof
\end{lemma}

As noticed by Runde (see \cite{Runde1}), there are very few
Banach algebras which are both dual and amenable.
For von Neumann algebras, which are the motivating example
of dual Banach algebras, there is a weaker notion of
amenablity, called Connes-amenability, which has a natural
generalisation to the case of dual Banach algebras.

\begin{definition}
Let $\mc A$ be a dual Banach algebra with predual $\mc A_*$.
Let $E$ be a Banach $\mc A$-bimodule.  Then $E'$ is a
\emph{w$^*$-Banach $\mc A$-bimodule} if, for each
$\mu\in E'$, the maps
\[ \mc A\rightarrow E',\quad a\mapsto
\begin{cases} a\cdot \mu, \\ \mu\cdot a \end{cases} \]
are $\sigma(\mc A,\mc A_*) - \sigma(E',E)$ continuous.

Then $(\mc A,\mc A_*)$ is Connes-amenable if, for each
w$^*$-Banach $\mc A$-bimodule $E'$, each derivation
$d:\mc A\rightarrow E'$, which is $\sigma(\mc A,\mc A_*) -
\sigma(E',E)$ continuous, is inner.
\noproof
\end{definition}

Given a Banach algebra $\mc A$, we define bilinear maps
$\mc A''\times\mc A'\rightarrow\mc A'$ and
$\mc A'\times\mc A''\rightarrow\mc A'$ by
\[ \ip{\Phi\cdot\mu}{a} = \ip{\Phi}{\mu\cdot a}
\quad \ip{\mu\cdot\Phi}{a} = \ip{\Phi}{a\cdot\mu}
\qquad (\Phi\in\mc A'', \mu\in\mc A', a\in\mc A). \]
We then define two bilinear maps $\aone,\atwo:\mc A''\times\mc A''
\rightarrow\mc A''$ by
\[ \ip{\Phi\aone\Psi}{\mu} = \ip{\Phi}{\Psi\cdot\mu}
\quad \ip{\Phi\atwo\Psi}{\mu} = \ip{\Psi}{\mu\cdot\Phi}
\qquad (\Phi,\Psi\in\mc A'', \mu\in\mc A'). \]
We can check that $\aone$ and $\atwo$ are actually algebra
products, called the \emph{first} and \emph{second Arens
products} respectively.  Then $\kappa_A:\mc A\rightarrow\mc A''$
is a homomorphism with respect to either Arens product.  When
$\aone = \atwo$, we say that $\mc A$ is \emph{Arens regular}.
In particular, when $\mc A$ is Arens regular, we may check that
$\mc A''$ is a dual Banach algebra with predual $\mc A'$.

\begin{theorem}\label{ca_facts}
Let $\mc A$ be an Arens regular Banach algebra.  When $\mc A$
is amenable, $\mc A''$ is Connes-amenable.  If $\tobidual{\mc A}$
is an ideal in $\mc A''$ and $\mc A''$ is Connes-amenable, then
$\mc A$ is amenable.

Let $\mc A$ be a C$^*$-algebra.  Then $\mc A$ is Arens regular,
and $\mc A''$ is Connes-amenable if and only if
$\mc A$ is amenable.
\end{theorem}
\begin{proof}
The first statements are \cite[Corollary~4.3]{Runde1} and
\cite[Theorem~4.4]{Runde1}.
The statement about C$^*$-algebras is detailed in
\cite[Chapter~6]{RundeBook}.
\end{proof}

Another class of Connes-amenable dual Banach algebras is
given by Runde in \cite{Runde4}, where it is shown that
$M(G)$, the measure algebra of a locally compact group $G$,
is amenable if and only if $G$ is amenable.

The organisation of this paper is as follows.  Firstly, we
study intrinsic characterisations of amenability, recalling
a result of Runde from \cite{Runde2}.  We then simplify these
conditions in the case of Arens regular Banach algebras.
We recall the notion of an \emph{injective} module, and quickly
note how Connes-amenability can be phrased in this language.
The final section of the paper then applies these ideas to
weighted semigroup algebras.  We finish with some open questions.

\section{Characterisations of amenability}

Let $E$ and $F$ be Banach spaces, and form the algebraic tensor
product $E\otimes F$.  We can norm $E\otimes F$ with the
\emph{projective tensor norm}, defined as
\[ \|u\|_\pi = \inf\Big\{ \sum_{k=1}^n \|x_k\|\|y_k\| :
u = \sum_{k=1}^n x_k \otimes y_k \Big\}
\qquad (u\in E\otimes F). \]
Then the completion of $(E\otimes F, \|\cdot\|_\pi)$ is
$E \proten F$, the \emph{projective tensor product} of $E$ and $F$.

Let $\mc A$ be a Banach algebra.  Then $\mc A\proten\mc A$ is a
Banach $\mc A$-bimodule for the module actions given by
\[ a\cdot (b\otimes c) = ab \otimes c,
\quad (b\otimes c) \cdot a = b\otimes ca
\qquad (a\in\mc A, b\otimes c\in\mc A\proten\mc A). \]
Define $\Delta_{\mc A} : \mc A \proten \mc A \rightarrow \mc A$ by
$\Delta_{\mc A}(a\otimes b)=ab$.  Then $\Delta_{\mc A}$ is an
$\mc A$-bimodule homomorphism.

\begin{theorem}\label{when_amen}
Let $\mc A$ be a Banach algebra.  Then the following are equivalent:
\begin{enumerate}
\item $\mc A$ is amenable;
\item $\mc A$ has a \emph{virtual diagonal}, which is a functional
   $M\in (\mc A\proten \mc A)''$ such that $a\cdot M =M\cdot a$ and
   $\Delta_{\mc A}''(M) \cdot a = \kappa_{\mc A}(a)$ for each $a\in\mc A$.
\end{enumerate}
\vspace{-2ex}\noproof
\end{theorem}

Runde introduced, in \cite{Runde2}, the following notion in
order to prove a version of the above theorem for Connes-amenability.

\begin{definition}
Let $\mc A$ be a dual Banach algebra with predual $\mc A_*$, and
let $E$ be a Banach $\mc A$-bimodule.  Then $x\in \sigma WC(E)$ if
and only if the maps $\mc A\rightarrow E$,
\[ a \mapsto \begin{cases} a\cdot x, \\ x\cdot a \end{cases} \]
are $\sigma(\mc A,\mc A_*) - \sigma(E,E')$ continuous.
\noproof
\end{definition}

It is clear that $\sigma WC(E)$ is a closed submodule of $E$.
The $\mc A$-bimodule homomorphism $\Delta_{\mc A}$ has adjoint
$\Delta'_{\mc A}:\mc A' \rightarrow (\mc A\proten\mc A)'$.
In \cite[Corollary~4.6]{Runde2} it is shown that
$\Delta'_{\mc A}(\mc A_*) \subseteq \sigma WC( (\mc A\proten\mc A)' )$.
Consequently, we can view $\Delta_{\mc A}'$ as a map
$\mc A_* \rightarrow \sigma WC( (\mc A\proten\mc A)' )$, and
hence view $\Delta''_{\mc A}$ as a map
$\sigma WC( (\mc A\proten\mc A)' )' \rightarrow \mc A_*'
= \mc A$, denoted by $\tilde \Delta_{\mc A}$.

\begin{theorem}\label{Runde_Thm}
Let $\mc A$ be a dual Banach algebra with predual $\mc A_*$.
Then the following are equivalent:
\begin{enumerate}
\item $\mc A$ is Connes-amenable;
\item $\mc A$ has a \emph{$\sigma WC$-virtual diagonal}, which
is $M\in \sigma WC( (\mc A\proten\mc A)' )'$ such that
$a\cdot M = M\cdot a$ and $a \tilde \Delta_{\mc A}(M) = a$ for
each $a\in\mc A$.
\end{enumerate}
\end{theorem}
\begin{proof}
This is \cite[Theorem~4.8]{Runde2}.
\end{proof}

In particular, we see that a Connes-amenable Banach algebra is
unital (which can of course be shown in an elementary fashion,
as in \cite[Proposition~4.1]{Runde1}).

\section{Connes-amenability for biduals of algebras}
\label{con_amen_bidual}

Recall Gantmacher's theorem, which states that a bounded linear
map $T:E\rightarrow F$ between Banach spaces $E$ and $F$ is
\emph{weakly-compact} if and only if $T''(E'') \subseteq \kappa_F(F)$.
We write $\mc W(E,F)$ for the collection of weakly-compact operators
in $\mc B(E,F)$.

\begin{lemma}\label{wsw_cty}
Let $E$ be a dual Banach space with predual $E_*$, let $F$ be a
Banach space, and let $T\in\mc B(E,F')$.  Then the following are
equivalent, and in particular each imply that $T$ is weakly-compact:
\begin{enumerate}
\item $T$ is $\sigma(E,E_*) - \sigma(F',F'')$ continuous;
\item $T'(F'') \subseteq \kappa_{E_*}(E_*)$;
\item there exists $S\in\mc W(F,E_*)$ such that $S' = T$.
\end{enumerate}
\end{lemma}
\begin{proof}
That (1) and (2) are equivalent is standard (compare with
Lemma~\ref{weak_star_cts}).

Suppose that (2) holds, so that we may define $S\in\mc B(F,E_*)$
by $\kappa_{E_*}\circ S = T'\circ\kappa_F$.  Then, for $x\in E$
and $y\in F$, we have
\[ \ip{x}{S(y)} = \ip{T'(\kappa_F(y))}{x}
= \ip{T(x)}{y}, \]
so that $S'=T$.  Then $S''(F'') = T'(F'') \subseteq
\kappa_{E_*}(E_*)$, so that $S$ is weakly-compact, by
Gantmacher's Theorem, so that (3) holds.

Conversely, if (3) holds, as $S$ is weakly-compact, we
have $\kappa_{E*}(E_*) \supseteq S''(F'') = T'(F'')$, so
that (2) holds.
\end{proof}

It is standard that for Banach spaces $E$ and $F$, we have
$(E\proten F)' = \mc B(F,E')$ with duality defined by
\[ \ip{T}{x\otimes y} = \ip{T(y)}{x}
\qquad (T\in\mc B(F,E'), x\otimes y\in E\proten F). \]
Then we see, for $a,b,c\in\mc A$ and $T\in(\mc A
\proten\mc A)' = \mc B(\mc A,\mc A')$, that
$\ip{a\cdot T}{b\otimes c} = \ip{T(ca)}{b}$ and that
$\ip{T\cdot a}{b\otimes c} = \ip{T(c)}{ab}
= \ip{T(c)\cdot a}{b}$ so that
\begin{equation}
(a\cdot T)(c) = T(ca), \quad (T\cdot a)(c) = T(c)\cdot a
\qquad (a,c\in\mc A, T:\mc A\rightarrow\mc A').
\label{eq:one} \end{equation}

Notice that we could also have defined $(E\proten F)'$ to be
$\mc B(E,F')$.  This would induce a different bimodule structure
on $\mc B(\mc A,\mc A')$, and we shall see in Section~\ref{Inj_predual}
that our chosen convention seems more natural for the task at hand.

\begin{proposition}\label{first_wap_prop}
Let $\mc A$ be a dual Banach algebra with predual $\mc A_*$.
For $T\in\mc B(\mc A,\mc A') = (\mc A\proten\mc A)'$, define
maps $\phi_r, \phi_l : \mc A\proten\mc A\rightarrow\mc A'$ by
\[ \phi_r(a\otimes b) = T'\kappa_{\mc A}(a) \cdot b,
\quad \phi_l(a\otimes b) = a \cdot T(b)
\qquad (a\otimes b\in\mc A\proten\mc A). \]
Then $T\in\sigma WC(\mc B(\mc A,\mc A'))$ if and only if
$\phi_r$ and $\phi_l$ are weakly-compact and have ranges
contained in $\tobidual{\mc A_*}$.
\end{proposition}
\begin{proof}
For $T\in\mc B(\mc A,\mc A') = (\mc A\proten\mc A)'$, define
$R_T,L_T : \mc A\rightarrow (\mc A\proten\mc A)'$ by
$R_T(a) = a\cdot T$ and $L_T = T\cdot a$, for $a\in\mc A$.
By definition, $T\in\sigma WC(\mc B(\mc A,\mc A'))$ if and only if
$R_T$ and $L_T$ are $\sigma(\mc A,\mc A_*)-%
\sigma(\mc B(\mc A,\mc A'),(\mc A\proten\mc A)'')$ continuous.
By Lemma~\ref{wsw_cty}, this is if and only if there
exist $\varphi_r, \varphi_l \in \mc W(\mc A\proten\mc A,\mc A_*)$
such that $\varphi_r' = R_T$ and $\varphi_l' = L_T$.

For $a\otimes b\in\mc A\proten\mc A$ and $c\in\mc A$, we see that
\begin{align*}
\ip{c}{\varphi_r(a\otimes b)} &= \ip{R_T(c)}{a\otimes b}
   = \ip{c\cdot T}{a\otimes b} = \ip{T(bc)}{a} \\
&= \ip{T'\kappa_{\mc A}(a)}{bc}
   = \ip{T'\kappa_{\mc A}(a)\cdot b}{c}
   = \ip{\phi_r(a\otimes b)}{c}, \\
\ip{c}{\varphi_l(a\otimes b)} &= \ip{L_T(c)}{a\otimes b}
   = \ip{T\cdot c}{a\otimes b} = \ip{T(b)}{ca} \\
&= \ip{a\cdot T(b)}{c} = \ip{\phi_l(a\otimes b)}{c}.
\end{align*}
Thus $\kappa_{\mc A_*}\circ\varphi_r = \phi_r$ and
$\kappa_{\mc A_*}\circ\varphi_l = \phi_l$.  Consequently,
we see that $T\in\sigma WC(\mc B(\mc A,\mc A'))$ if and only if
$\phi_r$ and $\phi_l$ are weakly-compact and take values in
$\tobidual{\mc A_*}$.
\end{proof}

The following definition is \cite[Definition~4.1]{Runde2}.

\begin{definition}
Let $\mc A$ be a Banach algebra and let $E$ be a Banach
$\mc A$-bimodule.  An element $x\in E$ is \emph{weakly almost
periodic} if the maps
\[ \mc A\rightarrow E,\quad
a \mapsto \begin{cases} a\cdot x, \\ x\cdot a \end{cases} \]
are weakly-compact.  The collection of weakly almost periodic
elements in $E$ is denoted by $\wap(E)$.
\noproof\end{definition}

\begin{lemma}\label{wap_to_maps}
Let $\mc A$ be a Banach algebra, and let $T\in\mc B(\mc A,\mc A')
= (\mc A\proten\mc A)'$. Let $\phi_r,\phi_l : \mc A\proten\mc A
\rightarrow\mc A'$ be as above.  Then
$T\in\wap( \mc B(\mc A,\mc A') )$ if and only if
$\phi_r$ and $\phi_l$ are weakly-compact.
\end{lemma}
\begin{proof}
Let $R_T,L_T:\mc A\rightarrow\mc B(\mc A,\mc A')$ be as in the
above proof.  By definition, $T\in\wap( \mc B(\mc A,\mc A') )$
if and only if $L_T$ and $R_T$ are weakly-compact.
We can verify that
\[ \phi_r'\circ\kappa_{\mc A} = R_T,\
\phi_l'\circ\kappa_{\mc A} = L_T,\
R_T'\circ\kappa_{\mc A\proten\mc A} = \phi_r,\
L_T'\circ\kappa_{\mc A\proten\mc A} = \phi_l, \]
which completes the proof.
\end{proof}

\begin{corollary}\label{unital_dual_sigma}
Let $\mc A$ be a unital, dual Banach algebra with predual $\mc A_*$,
and let $T\in\mc B(\mc A,\mc A') = (\mc A\proten\mc A)'$.  The
following are equivalent, and, in particular, each imply that $T$
is weakly-compact:
\begin{enumerate}
\item $T\in\sigma WC(\mc B(\mc A,\mc A'))$;
\item $T(\mc A) \subseteq
   \kappa_{\mc A_*}(\mc A_*)$, $T'(\kappa_{\mc A}(\mc A))
   \subseteq \kappa_{\mc A_*}(\mc A_*)$, and
   $T\in\sigma WC(\mc B(\mc A,\mc A'))$;
\item $T(\mc A) \subseteq
   \kappa_{\mc A_*}(\mc A_*)$, $T'(\kappa_{\mc A}(\mc A))
   \subseteq \kappa_{\mc A_*}(\mc A_*)$, and
   $T \in \wap( \mc B(\mc A,\mc A') )$.
\end{enumerate}
\end{corollary}
\begin{proof}
Let $e_{\mc A}$ be the unit of $\mc A$, so that for $a\in\mc A$,
we have $T(a) =	\phi_l(e_{\mc A} \otimes a)$ and
$T'\kappa_{\mc A}(a) = \phi_r(a \otimes e_{\mc A})$, which shows
that (1) implies (2); clearly (2) implies (1).

As $\mc A_*$ is an $\mc A$-bimodule, (2) and (3) are
equivalent by an application of Lemma~\ref{wap_to_maps}
and Proposition~\ref{first_wap_prop}.
\end{proof}

\begin{theorem}\label{When_Con_Amen}
Let $\mc A$ be a dual Banach algebra with predual $\mc A_*$.
Then $\mc A$ is Connes-amenable if and only if $\mc A$ is
unital and there exists
$M\in(\mc A\proten\mc A)''$ such that:
\begin{enumerate}
\item $\ip{M}{a\cdot T-T\cdot a}=0$ for $a\in\mc A$ and
   $T\in\sigma WC( \mc W(\mc A,\mc A') )$;
\item $\kappa_{\mc A_*}' \Delta''_{\mc A}(M) = e_{\mc A}$,
   where $e_{\mc A}$ is the unit of $\mc A$.
\end{enumerate}
\end{theorem}
\begin{proof}
As $\sigma WC( (\mc A\proten\mc A)' )'$ is a quotient of
$(\mc A\proten\mc A)''$, this is just a re-statement of
Theorem~\ref{Runde_Thm}.
\end{proof}

When $\mc A$ is an Arens regular Banach algebra, $\mc A''$
is a dual Banach algebra with canonical predual $\mc A'$.
In this case, we can make some significant simplifications
in the characterisation of when $\mc A''$ is Connes-amenable.

For a Banach algebra $\mc A$, we define the map
$\kappa_{\mc A} \otimes \kappa_{\mc A}: \mc A\proten \mc A
\rightarrow \mc A'' \proten \mc A''$ by
\[ ( \kappa_{\mc A} \otimes \kappa_{\mc A} ) (a\otimes b)
= \kappa_{\mc A}(a) \otimes \kappa_{\mc A}(b)
\qquad (a\otimes b\in \mc A\proten \mc A). \]
We turn $\mc A'' \proten \mc A''$ into a Banach $\mc A$-bimodule
in the canonical way.
Then $\kappa_{\mc A}\otimes\kappa_{\mc A}$ is an $\mc A$-bimodule
homomorphism.  The following is a simple verification.

\begin{lemma}\label{can_proj}
Let $\mc A$ be a Banach algebra.  The map
\[ \iota_{\mc A}:\mc B(\mc A,\mc A') \rightarrow
\mc B(\mc A'',\mc A'''); \ T\mapsto T'', \]
is an $\mc A$-bimodule homomorphism which is an
isometry onto its range.  Furthermore, we have that
$(\kappa_{\mc A} \otimes \kappa_{\mc A})' \circ \iota_{\mc A} =
I_{\mc B(\mc A,\mc A')}$.  Define $\rho_{\mc A}:
\mc A''\proten\mc A''\rightarrow (\mc A\proten\mc A)''$ by
\[ \ip{\rho_{\mc A}(\tau)}{T} = \ip{T''}{\tau}
\qquad (\tau\in \mc A''\proten \mc A'',
T\in\mc B(\mc A,\mc A')=(\mc A\proten \mc A)'). \]
Then $\rho_{\mc A}$ is a norm-decreasing $\mc A$-bimodule
homomorphism which satisfies $\rho_{\mc A} \circ
(\kappa_{\mc A}\otimes\kappa_{\mc A}) = \kappa_{\mc A\proten\mc A}$.
\noproof
\end{lemma}

For a Banach algebra $\mc A$, it is clear that $\mc W(\mc A,\mc A')$
is a sub-$\mc A$-bimodule of $\mc B(\mc A,\mc A')=(\mc A\proten\mc A)'$.

\begin{theorem}\label{arens_wap}
Let $\mc A$ be an Arens regular Banach algebra such that $\mc A''$
is unital, and let $T\in\mc B(\mc A'',\mc A''') =
(\mc A''\proten\mc A'')'$.  Then the following are equivalent:
\begin{enumerate}
\item $T \in \sigma WC(\mc B(\mc A'',\mc A'''))$, where we treat
   $\mc B(\mc A'',\mc A''')$ as an $\mc A''$-bimodule;
\item $T = S''$ for some $S\in \wap( \mc W(\mc A,\mc A') )$,
   where now we treat $\mc W(\mc A,\mc A')$ as an $\mc A$-bimodule.
\end{enumerate}
\end{theorem}
\begin{proof}
We apply Corollary~\ref{unital_dual_sigma} to $\mc A''$, so that
(1) is equivalent to $T$ being weakly-compact,
$T(\mc A'')\subseteq \kappa_{\mc A'}(\mc A')$,
$T'(\kappa_{\mc A''}(\mc A'')) \subseteq \kappa_{\mc A'}(\mc A')$,
and $T\in\wap( \mc B(\mc A'',\mc A''') )$.  Thus, if (1) holds,
then there exists $T_0 \in \mc W(\mc A'',\mc A')$ such that
$T = \kappa_{\mc A'} \circ T_0$, and there exists
$T_1 \in \mc W(\mc A'',\mc A')$ such that
$T'\circ\kappa_{\mc A''} = \kappa_{\mc A'} \circ T_1$.
Let $S = T_0\circ\kappa_{\mc A} \in \mc W(\mc A,\mc A')$.
Then, for $a\in\mc A$ and $\Psi\in\mc A''$, we have
\begin{align*}
\ip{S'(\Psi)}{a} &= \ip{\Psi}{T_0(\kappa_{\mc A}(a))}
= \ip{T(\kappa_{\mc A}(a))}{\Psi}
= \ip{T'(\kappa_{\mc A''}(\Psi))}{\kappa_{\mc A}(a)} \\
&= \ip{\kappa_{\mc A}(a)}{T_1(\Psi)}
= \ip{T_1(\Psi)}{a},
\end{align*}
so that $S' = T_1$.  Thus, for $\Phi,\Psi\in\mc A''$, we have
\begin{align*}
\ip{S''(\Phi)}{\Psi} &= \ip{\Phi}{T_1(\Psi)}
= \ip{T'(\kappa_{\mc A''}(\Psi))}{\Phi}
= \ip{T(\Phi)}{\Psi},
\end{align*}
so that $S''=T$.  We know that the maps $L_T,R_T:
\mc A''\rightarrow\mc B(\mc A'',\mc A''')$, defined by
$L_T(\Phi) = T\cdot\Phi$ and $R_T(\Phi) = \Phi\cdot T$ for
$\Phi\in\mc A''$, are weakly-compact.  Define $L_S,R_S:
\mc A\rightarrow\mc B(\mc A,\mc A')$ is an analogous manner,
using $S\in\mc W(\mc A,\mc A')$.  For $a\in\mc A$, $S\cdot a
\in\mc W(\mc A,\mc A')$, so for $\Psi\in\mc A''$ and $b\in\mc A$,
\[ \ip{(S\cdot a)'(\Psi)}{b} = \ip{\Psi}{(S\cdot a)(b)}
= \ip{\Psi}{S(b)\cdot a} = \ip{a\cdot\Psi}{S(b)}
= \ip{S'(a\cdot\Psi)}{b}. \]
Thus, for $a\in\mc A$ and
$\Phi,\Psi\in\mc A''$, we have that
\begin{align*}
\ip{ \iota_{\mc A}(L_S(a))(\Phi) }{\Psi} &=
\ip{ (S\cdot a)''(\Phi) }{\Psi} = \ip{ \Phi }{ S'(a\cdot\Psi) }
= \ip{ S''(\Phi) \cdot a }{\Psi},
\end{align*}
so that $\iota_{\mc A}(L_S(a))(\Phi) = S''(\Phi)\cdot a$,
and hence that $\iota_{\mc A}(L_S(a)) = S''\cdot a
= T\cdot a = T\cdot\kappa_{\mc A}(a) = L_T(\kappa_{\mc A}(a))$.
Thus we have that $L_S = (\kappa_{\mc A}\otimes\kappa_{\mc A})'
\circ R_T \circ \kappa_{\mc A}$, so that $L_S$
is weakly-compact.  A similar calculation shows that
$R_S$ is also weakly-compact, so that $S \in\wap(\mc W(\mc A,\mc A'))$.
This shows that (1) implies (2).

Conversely, if (2) holds, then $L_S$ and $R_S$ are weakly-compact.
As $S$ is weakly-compact, $T(\mc A'') = S''(\mc A'') \subseteq
\kappa_{\mc A'}(\mc A')$ and $T'(\kappa_{\mc A''}(\mc A''))
= S'''(\kappa_{\mc A''}(\mc A'')) = \kappa_{\mc A'}(S'(\mc A''))
\subseteq \kappa_{\mc A'}(\mc A')$, and $T$ is weakly-compact.
Thus, to show (1), we are required to show that $L_T$ and $R_T$
are weakly-compact.

For $a,b\in\mc A$ and $\Phi\in\mc A'$, we have
\[ \ip{(a\cdot S)'(\Phi)}{b} = \ip{\Phi}{S(ba)}
= \ip{a \cdot S'(\Phi)}{b}. \]
Then, for $\Phi,\Psi\in\mc A''$ and $a\in\mc A$,
we thus have
\begin{align*}
\ip{R_S'(\rho_{\mc A}(\Phi\otimes\Psi))}{a} &=
\ip{(a\cdot S)''}{\Phi\otimes\Psi} =
\ip{(a\cdot S)''(\Psi)}{\Phi} = \ip{\Psi}{a\cdot S'(\Phi)} \\
&= \ip{\Psi\cdot a}{S'(\Phi)} =
\ip{\Psi\aone\kappa_{\mc A}(a)}{S'(\Phi)}
= \ip{S'(\Phi)\cdot\Psi}{a}.
\end{align*}
Hence we see that $R_S'(\rho_{\mc A}(\Phi\otimes\Psi))
= S'(\Phi)\cdot\Psi$.  Let $U=R_S'\circ\rho_{\mc A}:
\mc A''\proten\mc A'' \rightarrow \mc A'$, so that as
$R_S$ is weakly-compact, so is $U$.  Then, for $\Phi,\Psi,
\Gamma\in\mc A''$, we have that
\begin{align*}
\ip{U'(\Gamma)}{\Phi\otimes\Psi} &=
\ip{\Gamma}{S'(\Phi)\cdot\Psi} = \ip{\Psi\atwo\Gamma}{S'(\Phi)}
= \ip{S''(\Psi\aone\Gamma)}{\Phi}
= \ip{(\Gamma\cdot S'')(\Psi)}{\Phi},
\end{align*}
so that $U'(\Gamma) = \Gamma\cdot T$, that is,
$U' = R_T$, so that $R_T$ is weakly-compact.
Similarly, we can show that
$L_T$ is weakly-compact, completing the proof.
\end{proof}

\begin{theorem}\label{Connes_Amen}
Let $\mc A$ be an Arens regular Banach algebra.  Then $\mc A''$
is Connes-amenable if and only if $\mc A''$ is unital and there
exists $M \in (\mc A \proten \mc A)''$ such that:
\begin{enumerate}
\item $\Delta_{\mc A}''(M) = e_{\mc A''}$, the unit of $\mc A''$;
\item $\ip{M}{a\cdot T-T\cdot a}=0$ for each $a\in\mc A$ and each
   $T\in \wap( \mc W(\mc A,\mc A') )$.
\end{enumerate}
\end{theorem}
\begin{proof}
By Theorem~\ref{When_Con_Amen}, we wish to show that the existence
of such an $M$ is equivalent to the existence of
$N\in(\mc A''\proten\mc A'')''$ such that:
\begin{description}
\item[(N1)] $\kappa_{\mc A'}'\Delta_{\mc A''}''(N) = e_{\mc A''}$;
\item[(N2)] $\ip{N}{\Phi\cdot S - S\cdot\Phi}=0$ for each $\Phi\in\mc A''$
   and each $S\in\sigma WC( \mc B(\mc A'',\mc A''') )$.
\end{description}

We can verify that $\iota_{\mc A} \circ \Delta_{\mc A}'
= \Delta'_{\mc A''} \circ \kappa_{\mc A'}$, so that (N1)
is equivalent to $\Delta_{\mc A}'' \iota'_{\mc A}(N) =
e_{\mc A''}$.
For $S\in\sigma WC( \mc B(\mc A'',\mc A''') )$, we know that
$S=T''$ for some $T \in \wap( \mc W(\mc A,\mc A') )$,
by Theorem~\ref{arens_wap}.  That is, the maps $\phi_r$ and
$\phi_l$, formed using $T$ as in Proposition~\ref{first_wap_prop},
are weakly-compact.
Then, for $\Phi\in\mc A''$, $\phi_r'(\Phi), \phi_l'(\Phi) \in
\mc B(\mc A,\mc A')$, and we can check that
\[ \phi_r'(\Phi)(a) = \kappa_{\mc A}'T''(a\cdot\Phi),
\quad \phi_l'(\Phi)(a) = T(a)\cdot\Phi
\qquad (a\in\mc A). \]
Then $\phi_r'(\Phi)', \phi_l'(\Phi)'\in\mc B(\mc A'',\mc A')$
are the maps
\[ \phi_r'(\Phi)'(\Psi) = \Phi\cdot T'(\Psi),
\quad \phi_l'(\Phi)'(\Psi) = T'(\Phi\aone\Psi)
\qquad (\Psi\in\mc A''), \]
where we remember that $T''(\mc A'') \subseteq
\kappa_{\mc A'}(\mc A')$. Consequently
$\phi_r'(\Phi)'', \phi_l'(\Phi)''\in\mc B(\mc A'',\mc A''')$
are given by
\[ \phi_r'(\Phi)''(\Psi) = T''(\Psi\aone\Phi),
\quad\phi_l'(\Phi)''(\Psi) = T''(\Psi) \cdot \Phi
\qquad (\Psi\in\mc A''), \]
where $\mc A'''$ is an $\mc A''$-bimodule, as $\mc A''$
is Arens regular.
That is, $\phi_r'(\Phi)'' = \Phi\cdot S$ and
$\phi_l'(\Phi)'' = S\cdot\Phi$.  Hence (N2) is equivalent to
\[ 0 = \ip{N}{\phi_r'(\Phi)'' - \phi_l'(\Phi)''}
= \ip{N}{\iota_{\mc A}( \phi_r'(\Phi) - \phi_l'(\Phi) )}
= \ip{\iota_{\mc A}'(N)}{\phi_r'(\Phi) - \phi_l'(\Phi)}, \]
for each $\Phi\in\mc A''$ and
$S\in\sigma WC( \mc B(\mc A'',\mc A''') )$.
That is, (N2) is equivalent to
\[ \phi_r''\iota_{\mc A}'(N) - \phi_l''\iota_{\mc A}'(N) = 0
\qquad ( S\in\sigma WC( \mc B(\mc A'',\mc A''') ) ).  \]

As $\phi_r$ and $\phi_l$ are weakly-compact, $\phi_r''$ and
$\phi_l''$ take values in $\tobidual{\mc A'}$, and so (N2) is
equivalent to
\[ 0 = \ip{\phi_r''\iota_{\mc A}'(N) - \phi_l''\iota_{\mc A}'(N)}
{\kappa_{\mc A}(a)} =
\ip{\iota_{\mc A}'(N)}{ \phi_r'(\kappa_{\mc A}(a)) -
\phi_l'(\kappa_{\mc A}(a)) }, \]
for each $a\in\mc A$ and each
$S\in\sigma WC( \mc B(\mc A'',\mc A''') )$.
However, $\phi_r'(\kappa_{\mc A}(a)) -
\phi_l'(\kappa_{\mc A}(a)) = a\cdot T - T\cdot a$, so that (N2)
is equivalent to
\[ 0 = \ip{\iota_{\mc A}'(N)}{a\cdot T - T \cdot a}
\qquad (a\in\mc A), \]
for each $T\in\mc W(\mc A,\mc A')$ such that $\phi_r$ and $\phi_l$
are weakly-compact.

Thus we have established that (N1) holds for $N$ if and only if
(1) holds for $M=\iota_{\mc A}'(N)$, and that (N2) holds for $N$
if and only if (2) holds for $M=\iota_{\mc A}'(N)$, completing
the proof.
\end{proof}

We immediately see that $\mc A$ amenable implies that
$\mc A''$ is Connes-amenable.  Furthermore, if $\mc A$ is
itself a dual Banach algebra, then
Corollary~\ref{unital_dual_sigma} shows that if $\mc A''$
is Connes-amenable, then $\mc A$ is Connes-amenable: notice
that if $e_{\mc A''}$ is the unit of $\mc A''$, then
\[ \ip{\kappa_{\mc A_*}'(e_{\mc A''})a}{\mu}
= \ip{e_{\mc A''} \cdot a}{\kappa_{\mc A_*}(\mu)}
= \ip{\kappa_{\mc A}(a)}{\kappa_{\mc A_*}(\mu)}
= \ip{a}{\mu}
\qquad (a\in\mc A, \mu\in\mc A_*), \]
so that $\kappa_{\mc A_*}'(e_{\mc A''})$ is the unit of $\mc A$.

\section{Injectivity of the predual module}
\label{Inj_predual}

Let $\mc A$ be a Banach algebra, and let $E$ and $F$ be
Banach left $\mc A$-modules.  We write ${_{\mc A}}\mc B(E,F)$ for
the closed subspace of $\mc B(E,F)$ consisting of left $\mc A$-module
homomorphisms, and similarly write $\mc B_{\mc A}(E,F)$ and
${_{\mc A}}\mc B_{\mc A}(E,F)$ for right $\mc A$-module and
$\mc A$-bimodule homomorphisms, respectively.  We say that
$T\in{_{\mc A}}\mc B(E,F)$ is \emph{admissible} if both the kernel
and image of $T$ are closed, complemented subspaces of, respectively,
$E$ and $F$.  If $T$ is injective, this is equivalent to the existence
of $S\in\mc B(F,E)$ such that $ST=I_E$.

\begin{definition}
Let $\mc A$ be a Banach algebra, and let $E$ be a Banach left
$\mc A$-module.  Then $E$ is \emph{injective} if, whenever $F$ and $G$
are Banach left $\mc A$-modules, $\theta\in{_{\mc A}}\mc B(F,G)$ is
injective and admissible, and $\sigma\in{_{\mc A}}\mc B(F,E)$,
there exists $\rho\in{_{\mc A}}\mc B(G,E)$ with $\rho\circ\theta = \sigma$.
\end{definition}

We say that $E$ is \emph{left-injective} when we wish to stress that
we are treating $E$ as a left module.  Similar definitions hold for
right modules and bimodules (written \emph{right-injective} and
\emph{bi-injective} where necessary).

Let $\mc A$ be a Banach algebra, let $E$ be a Banach left
$\mc A$-module, and turn $\mc B(\mc A,E)$ into a left $\mc A$-module
by setting
\[ (a\cdot T)(b) = T(ba) \qquad (a,b\in\mc A, T\in\mc B(\mc A,E)). \]
Then there is a canonical left $\mc A$-module homomorphism
$\iota:E\rightarrow \mc B(\mc A,E)$ given by
\[ \iota(x)(a) = a\cdot x\qquad (a\in\mc A,x\in E). \]
Notice that if $E$ is a closed submodule of $\mc A'$, then
$\mc B(\mc A,E)$ is a closed submodule of $(A\proten A)' =
\mc B(\mc A,\mc A')$, and $\iota$ is the restriction of
$\Delta'_{\mc A}:\mc A'\rightarrow \mc B(\mc A,\mc A')$ to $E$.

Similarly, we turn $\mc B(\mc A\proten\mc A,E)$ into a Banach
$\mc A$-bimodule by
\[ (a\cdot T)(b\otimes c) = T(ba\otimes c), \
(T\cdot a)(b\otimes c) = T(b\otimes ac) \quad
( a,b,c\in\mc A, T\in\mc B(\mc A\proten\mc A,E) ). \]
We then define (with an abuse of notation)
$\iota:E\rightarrow\mc B(\mc A\proten\mc A,E)$ by
\[ \iota(x)(a\otimes b) = a\cdot x \cdot b\qquad (x\in E,
a\otimes b\in\mc A\proten\mc A), \]
so that $\iota$ is an $\mc A$-bimodule homomorphism.

We can also turn $\mc B(\mc A,E)$ into a right $\mc A$-module
by reversing the above (in particular, we need to take the
other possible choice in Section~\ref{con_amen_bidual} leading
to different module actions as compared to those in (\ref{eq:one}).)

\begin{proposition}
Let $\mc A$ be a Banach algebra, and let $E$ be a \emph{faithful}
Banach left $\mc A$-module (that is, for each non-zero $x\in E$
there exists $a\in\mc A$ with $a\cdot x\not=0$).  Then $E$ is
injective if and only if there exists
$\phi\in{_{\mc A}}\mc B( \mc B(\mc A,E), E )$ such that
$\phi \circ \iota = I_E$.

Similarly, if $E$ is a left and right faithful Banach $\mc A$-bimodule
(that is, for each non-zero $x\in E$ there exists $a,b\in\mc A$
with $a\cdot x\not=0$ and $x\cdot b\not=0$).  Then $E$ is 
injective if and only if there exists $\phi\in
{_{\mc A}}\mc B_{\mc A}( \mc B(\mc A\proten\mc A,E), E)$ such that
$\phi\circ\iota = I_E$.
\end{proposition}
\begin{proof}
The first claim is \cite[Proposition~1.7]{DP}, and the second
claim is an obvious generalisation.
\end{proof}

Again, there exists a similar characterisation for right modules.

Let $\mc A$ be a dual Banach algebra with predual $\mc A_*$.
It is simple to show (see \cite{Runde2}) that if $\mc A_*$ is
bi-injective, then $\mc A$ is Connes-amenable.  Helemskii showed in
\cite{Hel2} that for a von Neumann algebra $\mc A$, the converse
is true.  However, Runde (see \cite{Runde2}) and Tabaldyev
(see \cite{Tab}) have shown that $M(G)$, the measure algebra of a
locally compact group $G$, while being a dual Banach algebra with
predual $C_0(G)$, has that $C_0(G)$ is a left-injective $M(G)$-module
only when $G$ is finite.  Of course, Runde (see \cite{Runde4}) has
shown that $M(G)$ is Connes-amenable if and only if $G$ is amenable.

Similarly, it is simple to show (using a virtual diagonal) that
if $\mc A$ is a Banach algebra with a bounded approximate identity,
then $\mc A$ is amenable if and only if $\mc A'$ is bi-injective.

Let $E$ and $F$ be Banach left $\mc A$-modules, and let
$\phi:E\rightarrow F$ be a left $\mc A$-module homomorphism
which is bounded below.  Then $\phi(E)$ is a closed submodule of
$F$, so that $F / \phi(E)$ is a Banach left $\mc A$-module.
Hence we have a \emph{short exact sequence}:
\[ \spreaddiagramrows{4ex} \spreaddiagramcolumns{6ex}
\xymatrix{ 0 \ar[r] &
E \ar@<1ex>[r]^{\phi} &
F \ar@{->>}[r] \ar@<1ex>@{-->}[l]^{P} &
F/\phi(E) \ar[r] & 0 }. \]
If there exists a bounded linear map $P:F\rightarrow E$ such
that $P\circ\phi = I_E$, then we say that the short exact
sequence is \emph{admissible}.  If, further, we may choose $P$
to be a left $\mc A$-module homomorphism, then the short
exact sequence is said to \emph{split}.  Similar definitions
hold for right modules and bimodules.

\begin{proposition}
Let $\mc A$ be a Banach algebra, let $E$ be a Banach left
$\mc A$-module, and consider the following admissible
short exact sequence:
\[ \spreaddiagramrows{4ex} \spreaddiagramcolumns{6ex}
\xymatrix{ 0 \ar[r] &
E \ar@<1ex>[r]^{\iota} &
\mc B(\mc A,E) \ar@{->>}[r] \ar@<1ex>@{-->}[l]^{P} &
\mc B(\mc A,E) / \iota(E) \ar[r] & 0 }. \]
Then $E$ is injective if and only if this short exact sequence
splits.
\end{proposition}
\begin{proof}
See, for example, \cite[Section~5.3]{RundeBook}.
\end{proof}

\begin{proposition}
Let $\mc A$ be a unital dual Banach algebra with predual $\mc A_*$,
and consider the following admissible short exact sequence of
$\mc A$-bimodules:
\begin{equation} \spreaddiagramcolumns{3ex}
\xymatrix{ 0 \ar[r] &
\mc A_* \ar@<1ex>[r]^(.35){\Delta_{\mc A}'} &
\sigma WC((\mc A\proten\mc A)') \ar@{->>}[r] \ar@<1ex>@{-->}[l]^(0.65){P} &
\sigma WC((\mc A\proten\mc A)')/\Delta_{\mc A}'(\mc A_*) \ar[r] & 0 }.
\label{con_amen_ses} \end{equation}
Then $\mc A$ is Connes-amenable if and only if
this short exact sequence splits.
\end{proposition}
\begin{proof}
Notice that $\Delta'_{\mc A}$ certainly maps $\mc A_*$ into
$\sigma WC((\mc A\proten\mc A)') = \sigma WC(\mc B(\mc A,\mc A'))$,
and that Corollary~\ref{unital_dual_sigma} shows that we can define
$P:\sigma WC(\mc B(\mc A,\mc A'))\rightarrow\mc A_*$ by $P(T) =
T(e_{\mc A})$ for $T\in \sigma WC(\mc B(\mc A,\mc A'))$.

Suppose that we can choose $P$ to be an $\mc A$-bimodule homomorphism.
Then let $M = P'(e_{\mc A})$, so that for $a\in\mc A$ and $T\in\sigma
WC(\mc B(\mc A,\mc A'))$,
\[ \ip{a\cdot M-M\cdot a}{T} = \ip{e_{\mc A}}{P(T\cdot a-a\cdot T)}
= \ip{a-a}{P(T)} = 0, \]
so that $a\cdot M - M\cdot a$.  Also $\Delta''_{\mc A}(M) =
(P\circ\Delta'_{\mc A})'(e_{\mc A}) = e_{\mc A}$, so that $M$ is
a $\sigma WC$-virtual diagonal, and hence $\mc A$ is Connes-amenable
by Runde's theorem.

Conversely, let $M$ be a $\sigma WC$-virtual diagonal and define
$P:\sigma WC(\mc B(\mc A,\mc A'))\rightarrow\mc A'$ by
\[ \ip{P(T)}{a} = \ip{M}{a\cdot T} \qquad (a\in\mc A,
T\in\sigma WC(\mc B(\mc A,\mc A')). \]
Let $(a_\alpha)$ be a bounded net in $\mc A$ which tends to $a\in\mc A$
in the $\sigma(\mc A,\mc A_*)$-topology.  By definition,
$a_\alpha\cdot T\rightarrow a\cdot T$ weakly, for each
$T\in\sigma WC(\mc B(\mc A,\mc A'))$, so that
$\ip{P(T)}{a_\alpha} \rightarrow \ip{P(T)}{a}$.  This implies that
$P$ maps into $\mc A_*$, as required.  Then, for $\mu\in\mc A_*$,
\[ \ip{a}{P\Delta_{\mc A}'(\mu)} = \ip{M}{a\cdot\Delta_{\mc A}'(\mu)}
= \ip{M}{\Delta_{\mc A}'(a\cdot\mu)} = \ip{e_{\mc A}}{a\cdot \mu}
= \ip{a}{\mu} \qquad (a\in\mc A), \]
so that $P\Delta_{\mc A}'=I_{\mc A_*}$.  Finally, we note that
\begin{align*}
\ip{P(a\cdot T\cdot b)}{c} &= \ip{M}{ca\cdot T\cdot b}
= \ip{b\cdot M}{ca\cdot T}= \ip{M\cdot b}{ca\cdot T} \\
&= \ip{P(T)}{bca} = \ip{a\cdot P(T) \cdot b}{c}
\qquad (a,b,c\in\mc A, T\in\sigma WC(\mc B(\mc A,\mc A'))),
\end{align*}
so that $P$ is an $\mc A$-bimodule homomorphism, as required.
\end{proof}

Let $\mc A$ be an Arens regular Banach algebra.
By reversing the argument Theorem~\ref{arens_wap}, we can show
that $\Delta'_{\mc A}:\mc A'\rightarrow\mc B(\mc A,\mc A')$
actually maps into $\wap(\mc W(\mc A,\mc A'))$.  Furthermore,
if $\mc A''$ is unital, then we may define $P:\wap(\mc W(\mc A,\mc A'))
\rightarrow \mc A'$ by
\[ \ip{P(T)}{a} = \ip{e_{\mc A''}}{P(a)}
\qquad (a\in\mc A, T\in\wap(\mc W(\mc A,\mc A'))). \]
Then we have that
\[ \ip{P\Delta'_{\mc A}(\mu)}{a} = \ip{e_{\mc A''}}{a\cdot\mu}
= \ip{\mu}{a} \qquad (a\in\mc A, \mu\in\mc A'). \]

\begin{proposition}
Let $\mc A$ be an Arens regular Banach algebra such that
$\mc A''$ is unital, and consider the following admissible short
exact sequence of $\mc A$-bimodules:
\begin{equation} \spreaddiagramcolumns{3ex}
\xymatrix{ 0 \ar[r] &
\mc A' \ar@<1ex>[r]^(0.3){\Delta_{\mc A}'} &
\wap(\mc W(\mc A,\mc A')) \ar@{->>}[r] \ar@<1ex>[l]^(0.7){P} &
\wap(\mc W(\mc A,\mc A')) / \Delta_{\mc A}'(\mc A') \ar[r] & 0 }.
\label{eq:two} \end{equation}
Then $\mc A''$ is Connes-amenable if and only if this short exact
sequence splits.
\end{proposition}
\begin{proof}
This follows in the same manner as the above proof, using
Theorem~\ref{Connes_Amen}.
\end{proof}

\section{Beurling algebras}

Let $S$ be a discrete semigroup (we can extend the following
definitions to locally compact semigroups, but for the questions we
are interested in, the results for non-discrete groups are
trivial).  A \emph{weight} on $S$ is a function $\omega: S
\rightarrow \mathbb R^{>0}$ such that
\[ \omega(st)\leq \omega(s)\omega(t) \qquad (s,t\in S). \]
Furthermore, if $S$ is unital with unit $u_S$, then we also
insist that $\omega(u_S)=1$.  This last condition is simply
a normalisation condition, as we can always set $\hat\omega(s)
= \sup\{ \omega(st) \omega(t)^{-1} : t\in S \}$ for each $s\in S$.
For $s,t\in S$, we have that $\omega(st) \leq \hat\omega(s)\omega(t)$,
so that
\[ \hat\omega(st) = \sup\{ \omega(str)\omega(r)^{-1} : r\in S \}
\leq \sup \{ \hat\omega(s)\omega(tr)\omega(r)^{-1} : r\in S \}
= \hat\omega(s) \hat\omega(t). \]
Clearly $\hat\omega(u_S)=1$ and $\hat\omega(s)\leq\omega(s)$
for each $s\in S$, while $\hat\omega(s) \geq \omega(s)\omega(u_S)^{-1}$,
so that $\hat\omega$ is equivalent to $\omega$.

We form the Banach space
\[ l^1(S,\omega) = \Big\{ (a_g)_{g\in S} \subseteq\mathbb C :
\|(a_g)\| := \sum_{g\in S} |a_g| \omega(g) < \infty \Big\}. \]
Then $l^1(S,\omega)$, with the convolution product, is a Banach
algebra, called a \emph{Beurling algebra}.
See \cite{CY} and \cite{DL} for further information on Beurling
algebras and, in particular, their second duals.

It will be more convenient for us to think of $l^1(S,\omega)$
as the Banach space $l^1(S)$ together with a weighted
algebra product.  Indeed, for $g\in S$, let $\delta_g\in l^1(S)$
be the standard unit vector basis element which is thought of as
a point-mass at $g$.  Then each $x\in l^1(S)$ can be written
uniquely as $x = \sum_{g\in S} x_g \delta_g$ for some family
$(x_g) \subseteq \mathbb C$ such that $\|x\| = \sum_{g\in S}
|x_g| <\infty$.  We then define
\[ \delta_g \star_\omega \delta_h = \delta_g \star \delta_h
= \delta_{gh} \Omega(g,h)  \qquad (g,h\in S), \]
where $\Omega(g,h) = \omega(gh)\omega(g)^{-1}\omega(h)^{-1}$,
and extend $\star$ to $l^1(S)$ by linearity and continuity.

For example, if $\omega$ and $\hat\omega$ are equivalent weights
on $S$, the define $\psi: l^1(S,\omega) \rightarrow l^1(S,\hat\omega)$
by $\psi(\delta_s) = \hat\omega(s) \omega(s)^{-1} \delta_s$.  As $\omega$
and $\hat\omega$ are equivalent, $\psi$ is an isomorphism of
Banach spaces.  Then
$\psi(\delta_s \star \delta_t) = \omega(st) \omega(s)^{-1} \omega(t)^{-1}
\hat\omega(st) \omega(st)^{-1} \delta_{st} =
\psi(\delta_s) \star \psi(\delta_t)$, so that $\psi$ is a homomorphism.

For a set $I$, we define the space $c_0(I)$ as
\[ c_0(I) = \Big\{ (a_i)_{i\in I} : \forall\, \epsilon>0,
|\{ i\in I : |a_i|\geq\epsilon \}|<\infty \Big\}, \]
where $| \cdot |$ is the cardinality of a set.  We equip $c_0(I)$ with
the supremum norm; then $c_0(I)' = l^1(I)$.  For $i\in I$, we
let $e_i\in c_0(I)$ be the point mass at $i$, that is,
$\ip{\delta_j}{e_i} = \delta_{i,j}$, the Kronecker delta,
for $\delta_j\in l^1(I)$.  Then $c_0(I)$ is the closed linear span
of $\{ e_i : i\in I \}$.  We let $l^\infty(I)$ be the Banach space
of all bounded families $(a_i)_{i\in I}$, with the supremum norm.
Then $l^1(I)'=l^\infty(I)$, we can treat $c_0(I)$ as a subspace of
$l^\infty(I)$, and the map $\kappa_{c_0(I)}: c_0(I)\rightarrow
l^\infty(I)$ is just the inclusion map.

For a semigroup $S$ and $s\in S$, we define maps $L_s,R_s:S
\rightarrow S$ by
\[ L_s(t) = st, \quad R_s(t) = ts \qquad (t\in S). \]
If, for each $s\in S$, $L_s$ and $R_s$ are finite-to-one maps,
then we say that $S$ is \emph{weakly cancellative}.  When $L_s$
and $R_s$ are injective for each $s\in S$, we say that $S$ is
\emph{cancellative}.
When $S$ is abelian and cancellative, a construction going back
to Grothendieck shows that $S$ is a sub-semigroup of some abelian
group.  However, this can fail to hold for non-abelian semigroups.

\begin{proposition}
Let $S$ be a weakly cancellative semigroup, let $\omega$ be a
weight on $S$, and let $\mc A = l^1(S,\omega)$.  Then
$c_0(S) \subseteq l^\infty(S) = \mc A'$ is a sub-$\mc A$-module
of $\mc A'$, so that $l^1(S,\omega)$ is a dual Banach algebra
with predual $c_0(S)$.
\end{proposition}
\begin{proof}
For $g,h\in S$ and $a=(a_s)_{s\in S} \in l^1(S,\omega)$, we have
\[ \ip{e_g \cdot \delta_h}{a} = \ip{e_g}{\delta_h \star a}
= \ip{e_g}{\sum_{s\in S} a_s \delta_{hs} \Omega(h,s)}
= \sum_{\{ s\in S : hs = g\}} a_s \Omega(h,s). \]
As $S$ is weakly cancellative, there exists at most finitely
many $s\in S$ such that $hs=g$, so that $e_g \cdot \delta_h$
is a member of $c_0(S)$.  Thus we see that $c_0(S)$ is a right
sub-$\mc A$-module of $\mc A'$.
The argument on the left follows in an analogous manner.
\end{proof}

Notice that the above result will hold for some semigroups $S$
which are not weakly cancellative, provided that the weight behaves
in a certain way.  However, it would appear that the later results
do not easily generalise to the non-weakly cancellative case.

Following \cite[Definition~2.2]{DL}, we have the following definition.

\begin{definition}
Let $I$ and $J$ be non-empty infinite sets, and let
$f:I\times J\rightarrow \mathbb C$ be a function.  Then
\emph{$f$ clusters on $I\times J$} if
\[ \lim_{n\rightarrow\infty} \lim_{m\rightarrow\infty} f(x_m,y_n)
= \lim_{m\rightarrow\infty} \lim_{n\rightarrow\infty} f(x_m,y_n), \]
whenever $(x_m)\subseteq I$ and $(y_n)\subseteq J$ are sequences
of distinct elements, and both iterated limits exist.

Furthermore, \emph{$f$ $0$-clusters on $I\times J$} if $f$ clusters
on $I\times J$, and the iterated limits are always $0$, when they exist.
\noproof
\end{definition}

From now on we shall exclude the trivial case when our (semi-)group
is finite.

\begin{theorem}\label{B_AR}
Let $S$ be a discrete, weakly cancellative semigroup,
and let $\omega$ be a weight on $S$.  Then the following are equivalent:
\begin{enumerate}
\item $l^1(S,\omega)$ is Arens regular;
\item for sequences of distinct elements $(g_j)$ and
   $(h_k)$ in $S$, we have
   \[ \lim_{j\rightarrow\infty} \lim_{k\rightarrow\infty}
   \Omega(g_j,h_k) = 0, \]
   whenever the iterated limit exists;
\item $\Omega$ $0$-clusters on $S\times S$.
\end{enumerate}
\end{theorem}
\begin{proof}
That (1) and (2) are equivalent for cancellative semigroups
is \cite[Theorem~1]{CY}.  Close examination of the proof shows
that this holds for weakly cancellative semigroups as well.
That (1) and (3) are equivalent follows by generalising the
proof of \cite[Theorem~7.11]{DL}, which is essentially an
application of Grothendieck's criterion for an operator to
be weakly-compact.  Alternatively, it follows easily that
(2) and (3) are equivalent by considering the \emph{opposite
semigroup} to $S$ where we reverse the product.
\end{proof}

In \cite{CY} it is also shown that if $G$ is a discrete,
uncountable group, then $l^1(G,\omega)$ is not Arens regular
for any weight $\omega$.  Furthermore, by \cite[Theorem~2]{CY},
if $G$ is a non-discrete locally compact group, then $L^1(G,\omega)$
is never Arens regular.

We shall consider both the Connes-amenability of $l^1(S,\omega)''$
and $l^1(S,\omega)$ (with respect to the canonical predual $c_0(S)$)
as, with reference to Corollary~\ref{unital_dual_sigma} and
Theorem~\ref{arens_wap}, the calculations should be similar.

\begin{proposition}\label{weak_comp_infty}
Let $I$ be a non-empty set, and let $X\subseteq l^\infty(I)$ be
a subset.  Then the following are equivalent:
\begin{enumerate}
\item $X$ is relatively weakly-compact;
\item $X$ is relatively sequentially weakly-compact;
\item the absolutely convex hull of $X$ is relatively weakly-compact;
\item if we define $f:I	\times X \rightarrow \mathbb C$ by
   $f(i,x) = \ip{x}{\delta_i}$ for $i\in I$ and $x\in X$, then
   $f$ clusters on $I\times X$;
\end{enumerate}
\end{proposition}
\begin{proof}
That (1) and (2) are equivalent is the Eberlien-Smulian theorem;
that (1) and (3) are equivalent is the Krein-Smulian theorem.
That (1) and (4) are equivalent is a result of Grothendieck,
detailed in, for example, \cite[Theorem~2.3]{DL}.
%
\end{proof}

It is standard that for non-empty sets $I$ and $J$, we have that
$l^1(I) \proten l^1(J) = l^1(I \times J)$, where, for $i\in I$
and $j\in J$, $\delta_i \otimes \delta_j\in l^1(I)\proten l^1(J)$
is identified with $\delta_{(i,j)}\in l^1(I\times J)$.  Thus we have
$( l^1(I) \proten l^1(J) )' = \mc B(l^1(I), l^\infty(J))
= l^1(I\times J)' = l^\infty(I\times J)$, where $T\in
\mc B(l^1(I), l^\infty(J))$ is identified with $(T_{(i,j)})
\in l^\infty(I\times J)$, where $T_{(i,j)} =
\ip{T(\delta_i)}{\delta_j}$.

\textbf{Is this paragraph used?}
Let $S$ be a countable, discrete, unital semigroup, and let
$\omega$ be a weight on $S$.  Then $l^1(S\times S)$ is a Banach
$l^1(S,\omega)$-bimodule, with module actions
\[ \delta_k \cdot \delta_{(g,h)} = \delta_{(kg,h)} \Omega(k,g)
\quad , \quad
\delta_{(g,h)} \cdot \delta_k = \delta_{(g,hk)} \Omega(h,k)
\qquad (g,h,k\in S). \]

For a non-empty set $I$, the unit ball of $l^1(I)$ is the closure of the
absolutely-convex hull of the set $\{ \delta_i : i\in I\}$, so that for
a Banach space $E$, by the Krein-Smulian theorem, a map
$T:l^1(I)\rightarrow E$ is weakly-compact if and only if the set
$\{ T(\delta_i) : i\in I \}$ is relatively weakly-compact in $E$.

\begin{proposition}\label{wap_c_zero}
Let $S$ be a weakly cancellative semigroup, let $\omega$ be a weight
on $S$, and let $\mc A = l^1(S,\omega)$.  Let $T\in\mc B(\mc A,\mc A')$
be such that $T(\mc A) \subseteq \kappa_{c_0(S)}(c_0(S))$ and
$T'(\kappa_{\mc A}(\mc A)) \subseteq \kappa_{c_0(S)}(c_0(S))$.
Then $T\in\mc W(\mc A,\mc A')$, and $T\in\wap(\mc W(\mc A,\mc A'))$
if and only if, for each sequence $(k_n)$ of distinct elements of $S$,
and each sequence $(g_m,h_m)$ of distinct elements of $S\times S$ such
that the repeated limits
\begin{gather}
\lim_n \lim_m \ip{T(\delta_{h_m})}{\delta_{k_ng_m}}, \
\lim_n \lim_m \Omega(k_n,g_m) \label{cond_one} \\
\lim_n \lim_m \ip{T(\delta_{h_mk_n})}{\delta_{g_m}}, \
\lim_n \lim_m \Omega(h_m,k_n) \label{cond_two}
\end{gather}
all exist, we have that at least one repeated limit in
each row is zero.
\end{proposition}
\begin{proof}
That $T$ is weakly-compact follows from Gantmacher's Theorem
(compare with Corollary~\ref{unital_dual_sigma}).
To show that $T\in\wap$, by Lemma~\ref{wap_to_maps}, we are required to
show that the maps $\phi_r$ and $\phi_l$ are weakly-compact.
We shall show that $\phi_l$ is weakly-compact if and only if one of
the repeated limits in the first line (\ref{cond_one}) is zero; the
proof that $\phi_r$ is related to (\ref{cond_two}) follows in a
similar way.  We have that
\[ \phi_l(\delta_{(g,h)}) = \phi_l(\delta_g \otimes \delta_h) = 
\delta_g \cdot T(\delta_h)   \qquad (g,h\in S). \]
By Proposition~\ref{weak_comp_infty}, $\phi_l$ is weakly-compact
if and only if the function
\[ S \times (S \times S) \rightarrow \mathbb C ; \
(k,(g,h)) \mapsto \ip{\delta_g\cdot T(\delta_h)}{\delta_k}
= \ip{T(\delta_h)}{\delta_{kg}} \Omega(k,g) \qquad (g,h,k\in S) \]
clusters on $S \times (S \times S)$.
As $T$ is weakly-compact, the function
\[ S \times S\rightarrow\mathbb C; \quad
(g,h) \mapsto \ip{T(\delta_g)}{\delta_h} \qquad (g,h\in S) \]
does cluster on $S\times S$.

Let $(k_n)$ be a sequence of distinct elements of $S$, and let
$(g_m,h_m)$ be a sequence of distinct elements of $S\times S$ such
that the iterated limits
\begin{equation}
\lim_n \lim_m \ip{T(\delta_{h_m})}{\delta_{k_ng_m}} \Omega(k_n,g_m)
, \quad
\lim_m \lim_n \ip{T(\delta_{h_m})}{\delta_{k_ng_m}} \Omega(k_n,g_m)
\label{eq:five} \end{equation}
exist.  We now investigate when these iterated limits are equal.

Suppose firstly that, by moving to a subsequence if necessary,
we have that $g_m = g$ for all $m$.  Further, by moving to a subsequence
if necessary, we may suppose that $\lim_n \Omega(k_n,g) = \alpha$, say,
and that $(k_ng)$ is a sequence of distinct elements (as $S$ is
weakly cancellative).  Then
\begin{align*}
\lim_n \lim_m \ip{T(\delta_{h_m})}{\delta_{k_ng_m}} & \Omega(k_n,g_m)
= \lim_n \Omega(k_n,g) \lim_m \ip{T(\delta_{h_m})}{\delta_{k_ng}} \\
&= \alpha \lim_n \lim_m \ip{T(\delta_{h_m})}{\delta_{k_ng}}
= \alpha \lim_m \lim_n \ip{T(\delta_{h_m})}{\delta_{k_ng}} \\
&= \lim_m \lim_n \ip{T(\delta_{h_m})}{\delta_{k_ng_m}} \Omega(k_n,g_m),
\end{align*}
where we can swap the order of taking limits, as $T$ is weakly-compact.

Alternatively, if we cannot move to a subsequence such that $(g_m)$ is
constant, then we may move to subsequence such that $(g_m)$ is a
sequence of distinct elements, and such that the iterated limits
\begin{gather*}
\lim_m \lim_n \Omega(k_n,g_m),\quad
\lim_n \lim_m \Omega(k_n,g_m), \\
\lim_m \lim_n \ip{T(\delta_{h_m})}{\delta_{k_ng_m}},\quad
\lim_n \lim_m \ip{T(\delta_{h_m})}{\delta_{k_ng_m}}
\end{gather*}
all exists.  As $T(\mc A)\subseteq\kappa_{c_0(S)}(c_0(S))$,
we have that
\[ \{ g\in S : |\ip{T(\delta_h)}{\delta_g}|\geq\epsilon \}
\text{ is finite } \qquad (\epsilon>0, h\in S). \]
Consequently, and using the fact that $S$ is weakly cancellative,
we see that
\[ \lim_n \ip{T(\delta_{h_m})}{\delta_{k_ng_m}} = 0 \]
for each $m$.  Hence the iterated limits in (\ref{eq:five})
are equal if and only if we have that at least one repeated
limit in (\ref{cond_one}) is zero.
\end{proof}

\begin{proposition}\label{auto_semi_wap}
Let $S$ be a discrete, unital, weakly cancellative semigroup, and let
$\omega$ be a weight on $S$ such that $\mc A = l^1(S,\omega)$ is
Arens regular.  Then $\wap(\mc W(\mc A,\mc A')) = \mc W(\mc A,\mc A')$.
\end{proposition}
\begin{proof}
Let $T\in\mc W(\mc A,\mc A')$.
We can follow the above proof through until the point at which we
use the fact that $T(\mc A)\subseteq \kappa_{c_0(S)}(c_0(S))$.
However, as $l^1(S,\omega)$ is Arens regular, by Theorem~\ref{B_AR},
we have that
\[ \lim_m \lim_n \Omega(k_n,g_m) =
\lim_n \lim_m \Omega(k_n,g_m) = 0, \]
so that the iterated limits in (\ref{eq:five}) must be $0$,
implying that $\phi_l$ is weakly-compact.
In a similar manner, $\phi_r$ is weakly-compact.
\end{proof}

\begin{theorem}\label{when_loneg_cam}
Let $S$ be a discrete weakly cancellative semigroup,
and let $\omega$ be a weight on $S$ such that $\mc A = l^1(S,\omega)$
is Arens regular and $\mc A''$ is unital with unit $e_{\mc A''}$.
Then $\mc A''$ is Connes-amenable if and only if
there exists $M\in (\mc A\proten\mc A)'' = l^\infty(S\times S)'$
such that:
\begin{enumerate}
\item $\ip{M}{(f_{gh}\Omega(g,h))_{(g,h)\in S\times S}} = \ip{e_{\mc A''}}{f}$
   for each bounded family $(f_g)_{g\in S}$;
\item $\ip{M}{( f(hk,g)\Omega(h,k) - f(h,kg)\Omega(k,g) )_
   {(g,h)\in S\times S}} = 0$ for each $k\in S$, and each bounded
   function $f:S\times S\rightarrow\mathbb C$ which clusters
   on $S\times S$.
\end{enumerate}
\end{theorem}
\begin{proof}
We use Theorem~\ref{Connes_Amen} and Proposition~\ref{auto_semi_wap}.
For $f=(f_g)_{g\in S} \in l^\infty(S)$, we have
\[ \ip{ \Delta'_{\mc A}(f) }{\delta_g\otimes \delta_h}
= \ip{f}{\delta_{gh}}\Omega(g,h) \qquad (g,h\in S), \]
so that $\Delta'_{\mc A}(f) = ( \ip{f}{\delta_{gh}}\Omega(g,h) )_{%
(g,h)\in S\times S} \in l^\infty(S\times S)$.  As $f\in l^\infty(S)$
was arbitrary, we have condition (1).

For $T\in\mc B(\mc A,\mc A')$, we treat $T$ as being a member
of $l^\infty(S\times S)$.  Then $T$ is weakly-compact if and
only if the family $(\ip{T(\delta_g)}{\delta_h})_{(g,h)\in S\times S}$
clusters on $S\times S$.  For $k\in S$, we have
\[ \ip{ \delta_k\cdot T - T\cdot \delta_k }{\delta_g\otimes \delta_h}
= \ip{T(\delta_{hk})}{\delta_g} \Omega(h,k)
- \ip{T(\delta_h)}{\delta_{kg}} \Omega(k,g). \]
Thus we have condition (2).
\end{proof}

Notice that if $S$ is unital with unit $u_S$, then the unit of
$\mc A$ (and hence $\mc A''$) is $\delta_{u_S}$.  In this case,
condition (1) reduces to $\ip{M}{(f_{gh}\Omega(g,h))_{(g,h)\in S\times S}}
= f_{u_S}$.

\begin{theorem}\label{when_loneg_am}
Let $S$ be a discrete unital semigroup, let $\omega$ be a weight
on $S$, and let $\mc A=l^1(S,\omega)$.
Then $\mc A$ is amenable if and only if there exists
$M\in (\mc A\proten\mc A)'' = l^\infty(S\times S)'$ such that:
\begin{enumerate}
\item $\ip{M}{(f_{gh}\Omega(g,h))_{(g,h)\in S\times S}} = f_{u_S}$,
   where $u_S\in S$ is the unit of $S$,
   for each bounded family $(f_g)_{g\in S}$;
\item $\ip{M}{( f(hk,g)\Omega(h,k) - f(h,kg)\Omega(k,g) )_
   {(g,h)\in S\times S}} = 0$ for each $k\in S$, and each
   bounded function $f:S\times S\rightarrow\mathbb C$.
\end{enumerate}
\end{theorem}
\begin{proof}
This follows from Theorem~\ref{when_amen} in the same way that
the above follows from Theorem~\ref{Connes_Amen}.
\end{proof}

Notice that condition (2) of Theorem~\ref{when_loneg_am}
is strictly stronger than condition (2) of
Theorem~\ref{when_loneg_cam}.

\begin{theorem}\label{C_amen_predual}
Let $S$ be a discrete, weakly cancellative semigroup, let $\omega$
be a weight on $S$, and let $\mc A = l^1(S,\omega)$ be unital with
unit $e_{\mc A}$.  Then $\mc A$ is Connes-amenable, with respect to
the predual $c_0(S)$, if and only if there exists
$M\in (\mc A\proten\mc A)'' = l^\infty(S\times S)'$ such that:
\begin{enumerate}
\item $\ip{M}{(f_{gh}\Omega(g,h))_{(g,h)\in S\times S}} =
   \ip{e_{\mc A}}{f}$ for each family $(f_g)_{g\in S} \in c_0(S)$;
\item $\ip{M}{( f(hk,g)\Omega(h,k) - f(h,kg)\Omega(k,g) )_
   {(g,h)\in S\times S}} = 0$ for each $k\in S$, and each
   bounded function $f:S\times S\rightarrow\mathbb C$ which satisfies
   the conclusions of Proposition~\ref{wap_c_zero}.
\end{enumerate}
\end{theorem}
\begin{proof}
We now use Theorem~\ref{When_Con_Amen}.  By $f$ satisfying
the conclusions of Proposition~\ref{wap_c_zero}, we identify
$f:S\times S\rightarrow\mathbb C$ with $T\in\mc B(\mc A,\mc A')$
by $\ip{T(\delta_g)}{\delta_h} = f(g,h)$, for $g,h\in S$.
\end{proof}

We shall now establish when $l^1(S,\omega)$ and $l^1(S,\omega)''$
are Connes-amenable.  For a discrete group $G$, a weight $\omega$
on $G$ and $h\in G$, define $J_h\in\mc B(l^\infty(G))$ by
\[ J_h(f) = \big( f_{hg} \Omega(h,g)\omega(h) \Omega(g^{-1},h^{-1})
\omega(h^{-1}) \big)_{g\in G} \qquad (f=(f_g)_{g\in G}\in l^\infty(G)). \]
Notice then that, for $f\in l^\infty(G)$, we have
\[ \| J_h(f) \| = \sup_g |f_{hg}| \omega(hg)\omega(g)^{-1}
   \omega(g^{-1}h^{-1})\omega(g^{-1})^{-1}
\leq \|f\| \omega(h) \omega(h^{-1}), \]
so that $J_h$ is bounded.

\begin{definition}\label{group_omega_amen}
Let $G$ be a discrete group, and let $\omega$ be a weight on $G$.
We say that $G$ is \emph{$\omega$-amenable} if there exists
$N\in l^\infty(G)'$ such that:
\begin{enumerate}
\item $\ip{N}{(\Omega(g,g^{-1}))_{g\in G}}=1$, where
   $\Omega$ is defined by $\omega$, and hence $(\Omega(g,g^{-1}))_{g\in G}$
   is a bounded family forming an element of $l^\infty(G)$;
\item $J_h'(N)=N$ for each $h\in G$.\vspace{-4ex}
\end{enumerate}
\noproof
\end{definition}

Notice that if $\omega$ is identically $1$, then this condition
reduces to the usual notion of a group being amenable (we usually
require that $N$ is a \emph{mean}, in that $N$ is a positive
functional on $l^\infty(G)$, but by forming real and imaginary
parts, and then positive and negative parts, we can easily generate
a non-zero scalar multiple of a mean from a functional $N$ satisfying
the definition above).

\begin{theorem}
Let $G$ be a discrete group, let $\omega$ be
a weight on $G$, and let $\mc A=l^1(G,\omega)$.  Then
the following are equivalent:
\begin{enumerate}
\item $\mc A$ is Connes-amenable, with respect to the predual $c_0(G)$;
\item $\mc A$ is amenable;
\item $G$ is $\omega$-amenable.
\end{enumerate}
Furthermore, if $\mc A$ is Arens regular, then these conditions
are equivalent to $\mc A''$ being Connes-amenable.
\end{theorem}
\begin{proof}
It is clear that (2) implies (1).  When $\mc A$ is Arens regular,
(2) implies that $\mc A''$ is Connes-amenable, and $\mc A''$
Connes-amenable implies (1).  We shall thus show that (1) implies
(3), and that (3) implies (2).

If (1) holds, then let $M\in l^\infty(G\times G)'$
be given as in Theorem~\ref{C_amen_predual}.  Define $\phi:l^\infty(G)
\rightarrow l^\infty(G\times G)$ by
\[ \ip{\phi(f)}{\delta_{(g,h)}} = \begin{cases}
f_g & : g = h^{-1}, \\ 0 &: g\not=h^{-1}, \end{cases}
\quad\qquad ( f=(f_g)_{g\in G} \in l^\infty(G) ). \]
Let $N = \phi'(M) \in l^\infty(G)'$.  Then we have
\[ \phi( (\Omega(g,g^{-1}))_{g\in G} ) = 
( \delta_{h,g^{-1}} \Omega(g,h) )_{(g,h)\in G\times G}
= ( \delta_{gh,e_G} \Omega(g,h) )_{(g,h)\in G\times G}, \]
where $\delta$ is the Kronecker delta, so that
\[ \ip{N}{(\Omega(g,g^{-1}))_{g\in G}} = \delta_{e_G,e_G} = 1, \]
by condition (1) on $M$ from Theorem~\ref{C_amen_predual};
clearly $(\delta_{e_G,g})_{g\in G} \in c_0(G)$.

Fix $k\in G$ and $f\in l^\infty(G)$.  Define $F:G\times G\rightarrow
\mathbb C$ by
\[ F(h,g) = \delta_{gh,k} f_g \omega(k) \omega(hk^{-1})
\omega(h)^{-1}.   \qquad (g,h\in G). \]
Then we have $|F(h,g)| \leq |f_g| |\omega(k)| |\omega(hk^{-1})|
|\omega(h)|^{-1} \leq \|f\|_\infty |\omega(k)| |\omega(k^{-1})|$,
so that $F$ is bounded.  Let $T:\mc A\rightarrow\mc A'$ be the operator
associated with $F$.  For $g,h\in G$, we have that $F(h,g)\not=0$
only when $gh=k$, so that $T(\mc A)\subseteq c_0(S)$ and
$T'(\kappa_{\mc A}(\mc A)) \subseteq c_0(S)$.
Furthermore, if $(k_n)$ is a sequence of distinct elements in $G$,
and $(g_m,h_m)$ is a sequence of distinct elements in $G\times G$, then
$\lim_n \lim_m F(h_m,k_ng_m) = 0$.  This follows, as for $n_0$ fixed,
$k_{n_0}g_mh_m = k$ only if $g_mh_m = k_{n_0}^{-1} k$, so if this holds
for all sufficiently large $m$, we have that $k_ng_mh_m\not=k$ for
sufficiently large $m$ and $n\not=n_0$.  Similarly, 
$\lim_n \lim_m F(h_mk_n,g_m) = 0$, so that $F$ satisfies the conditions
of Proposition~\ref{wap_c_zero}.

Notice that
\[ \ip{ \phi(J_k(f)) }{ \delta_{(g,h)} } =
\delta_{gh,e_G} \ip{ J_k(f) }{ \delta_g }
= \delta_{gh,e_G} f_{kg} \omega(kg)\omega(g)^{-1}
\omega(g^{-1}k^{-1})\omega(g^{-1})^{-1}. \]
Thus we have
\begin{align*}
F(hk,g)&\Omega(h,k) - F(h,kg)\Omega(k,g) \\
&= \delta_{ghk,k} f_g \omega(k) \omega(hkk^{-1})
\omega(hk)^{-1} \Omega(h,k) -
\delta_{kgh,k} f_{kg} \omega(k) \omega(hk^{-1})
\omega(h)^{-1} \Omega(k,g) \\
&= \delta_{gh,e_G} f_g -    \delta_{gh,e_G} f_{kg}
   \omega(hk^{-1}) \omega(h)^{-1} \omega(kg) \omega(g)^{-1} \\
&= \ip{\phi(f) - \phi(J_k(f)) }{ \delta_{(g,h)} }.
\end{align*}
So, by condition (2) from Theorem~\ref{C_amen_predual}, we have that
\[ \ip{N}{ f - J_k(f) } = 0, \]
which, as $f$ was arbitrary, shows that $N = J_k'(N)$, as required.

Now suppose that $G$ is $\omega$-amenable.  We shall show that
$\mc A$ is amenable, which completes the proof.
Define $\psi : l^\infty(G\times G)\rightarrow l^\infty(G)$ by
\[ \ip{ \psi(F) }{ \delta_g } = F(g,g^{-1})
\qquad ( F\in l^\infty(G\times G), g\in G). \]
Let $N\in l^\infty(G)'$ be as in Definition~\ref{group_omega_amen},
and let $M =\psi'(N)$.  Then let $(f_g)_{g\in G}$ be a bounded
family in $\mathbb C$, so that
\[ \ip{M}{(f_{gh}\Omega(g,h))_{(g,h)\in G\times G}} =
\ip{N}{(f_{e_G}\Omega(g,g^{-1}))_{g\in G}} = f_{e_G}, \]
verifying condition (1) of Theorem~\ref{when_loneg_am} for $M$.

Let $f:G\times G\rightarrow\mathbb C$ be a bounded function,
and let $k\in G$.  Then
\begin{align*}
\psi\big( (f(hk,g) & \Omega(h,k) - f(h,kg)\Omega(k,g))_{(g,h)\in G\times G} \big) \\
&= \big( f(g^{-1}k,g)\Omega(g^{-1},k) - f(g^{-1},kg)\Omega(k,g) \big)_{g\in G}.
\end{align*}
Define $F:G\times G\rightarrow\mathbb C$ by
\[ F(g,h) = f(hk,g) \Omega(h,k) \qquad (g,h\in G), \]
so that $F$ is bounded.  For $g\in G$, we have that
\begin{align*}
&\ip{ \psi(F)-J_k(\psi(F)) }{ \delta_g } \\ &=
f(g^{-1}k,g) \Omega(g^{-1},k)
- f((kg)^{-1}k,kg) \Omega((kg)^{-1},k) \omega(kg) \omega(g)^{-1}
   \omega(g^{-1}k^{-1}) \omega(g^{-1})^{-1} \\
&= f(g^{-1}k,g) \Omega(g^{-1},k)
- f(g^{-1},kg) \omega(k)^{-1} \omega(kg) \omega(g)^{-1} \\
&= f(g^{-1}k,g) \Omega(g^{-1},k)
- f(g^{-1},kg) \Omega(k,g).
\end{align*}
Consequently, using condition (2) of Definition~\ref{group_omega_amen},
we have established condition (2) of Theorem~\ref{when_loneg_am}
for $M$.  This shows that $l^1(G,\omega)$ is amenable.
\end{proof}

\begin{example}
If $S$ is a semigroup which is not cancellative, then it is possible
for $l^1(S)$ to be unital while $S$ is not.  For example, let
$S$ be $(\mathbb N,\max)$ (where $\mathbb N=\{1,2,3,\ldots\}$ say)
with adjoined idempotents $u$ and $v$ such that $uv=vu=1$ and
$un = nu = vn = nv = n$ for $n\in\mathbb N$.  Then $S$ is a
weakly cancellative, commutative semigroup without a unit, but
$e = \delta_u+\delta_v-\delta_1$ is easily seen to be a unit for
$l^1(S)$.  Indeed, $S$ is seen to be a finite semilattice of groups,
so by the result of \cite{Gron2}, $l^1(S)$ is amenable.
\end{example}

In \cite[Theorem~2.3]{Gron1} it is shown that if $l^1(S,\omega)$
is amenable for a cancellative, unital semigroup $S$ and some
weight $\omega$, then $S$ is actually a group.  We shall now show
that this holds for Connes-amenability as well.

For a cancellative, unital semigroup $S$, with unit $u_S$, if $g\in S$
is invertible, then $g$ has a unique inverse, denoted by $g^{-1}$.
Furthermore, if $g$ has a left inverse, say $hg=u_S$, then $ghg =
g = u_Sg$ so that $gh=u_S$; similarly, if $gh=u_S$ then $hg=u_S$.

\begin{theorem}
Let $S$ be a weakly cancellative semigroup, let $\omega$ be a weight
on $S$, and let $\mc A=l^1(S,\omega)$.  Suppose that $\mc A$ is
Connes-amenable with respect to the predual $c_0(S)$.
If $S$ is cancellative or unital, then $S$ is a group.
\end{theorem}
\begin{proof}
As $\mc A$ is Connes-amenable, let $M\in (\mc A\proten\mc A)''$
be as in Theorem~\ref{C_amen_predual}.  Then $\mc A$ is unital,
with unit $e_{\mc A} = (a_s)_{s\in S} \in l^1(S,\omega)$ say.
For now, we shall not assume that $e_{\mc A}$ has norm one, as
the standard renorming to ensure this will not (a priori)
necessarily yield an $l^1(S,\hat\omega)$ algebra for some
weight $\hat\omega$.  Suppose that $S$ is cancellative.
Fix $h\in S$, so that
\[ \sum_{s\in S} a_s \delta_{sh} \Omega(s,h)
= e_{\mc A}\star \delta_h =
\delta_h = \delta_h \star e_{\mc A} =
\sum_{s\in S} a_s \delta_{hs} \Omega(h,s). \]
In particular, for each $h\in S$ there is a unique $u_h\in S$
such that $h u_h = h$ (so that $h u_h h = h^2$ implying that
$u_h h=h$), and we have that $a_{u_h} \omega(u_h)^{-1} =1$.
We also see that $a_s=0$ for each $s\in S$ such that
$sh\not=h$, that is, $s\not=u_h$.  However, $h$ was arbitrary,
so that $S$ is unital with unit $u_S$, and $e_{\mc A} =
\omega(u_S) \delta_{u_S}$, where we can now assume that
$\omega(u_S)=1$ by a renorming.

Now suppose that $S$ is a unital, weakly cancellative semigroup,
so that the unit of $\mc A$ is $\delta_{u_S}$.
Suppose that $s\in S$ has no right inverse.  Define
$F:S\times S\rightarrow\mathbb C$ by
\[ F(h,sg)=0, \quad F(hs,g) = \begin{cases}
\Omega(g,hs) &: gh=u_S, \\ 0 &: \text{otherwise.} \end{cases}
\qquad (g,h\in S). \]
To show that this is well-defined, suppose that for
$g,h,j,k\in S$, we have that $h=js$, $sg=k$ and $kj=u_S$.
Then $s(gj) = kj = u_S$, so that $s$ has a right inverse,
a contradiction.  Then $F$ is bounded, so let $T:\mc A\rightarrow
\mc A'$ be the operator associated with $F$.  Then $F(a,b)\not=0$
only when $ba=s$, so as $S$ is weakly cancellative, we see that
$T(\mc A)\subseteq c_0(S)$ and $T'(\kappa_{\mc A}(\mc A))
\subseteq c_0(S)$. 

Suppose that for sequences of distinct elements $(k_n)\subseteq S$
and $(g_m,h_m)\subseteq S\times S$, we have that
\[ \lim_n \lim_m \ip{T(\delta_{h_m})}{\delta_{k_ng_m}}
= \lim_n \lim_m F(h_m,k_ng_m) \not=0. \]
Then, for some $N>0$ and $\epsilon>0$, for each $n\geq N$,
$\lim_m F(h_m,k_ng_m)\geq \epsilon$.  Hence, for $n\geq N$, there
exists $M_n>0$ such that if $m\geq M_n$, then $k_ng_mh_m = s$
(as otherwise $F(h_m,k_ng_m)=0$).  This, however, contradicts
$S$ being weakly cancellative.  Similarly, if
$\lim_n \lim_m \ip{T(\delta_{h_mk_n})}{\delta_{g_m}}\not=0$,
then we need $g_mh_mk_n=s$ for all $n,m$ sufficiently large,
which is a contradiction.  Thus $T$ satisfies all the conditions
of Proposition~\ref{wap_c_zero}.

Then, for $g,h\in S$, if $gh=u_S$, we have that $\Omega(h,s)\Omega(g,hs)
= \omega(h)^{-1} \omega(g)^{-1} = \Omega(g,h)$, so that
\[ F(hs,g)\Omega(h,s) - F(h,sg)\Omega(s,g)
= \begin{cases} \Omega(g,h) &: gh=u_S, \\ 0 &: \text{otherwise}.
\end{cases} \]
Hence condition (2) of Theorem~\ref{C_amen_predual} implies
that $\ip{M}{(\delta_{gh,u_S} \Omega(g,h))_{(g,h)\in S\times S}}=0$,
which contradicts condition (1) of this theorem.  Hence every
element of $S$ has a right inverse.

By symmetry (or by repeating the argument on the left) we see that
every element of $S$ has a left inverse, and that hence $S$ must
be a group.
\end{proof}

We hence have the following theorem, which shows that weighted
semigroup algebras behave like C$^*$-algebras with regards to
Connes-amenability.

\begin{theorem}\label{main_thm}
Let $S$ be a discrete cancellative semigroup, and let
$\omega$ be a weight on $S$.  The following are equivalent:
\begin{enumerate}
\item $l^1(S,\omega)$ is amenable;
\item $l^1(S,\omega)$ is Connes-amenable, with respect to
   the predual $c_0(S)$;
\end{enumerate}
If $l^1(S,\omega)$ is Arens regular, then
these conditions are equivalent to $l^1(S,\omega)''$ being
Connes-amenable.
These equivalent conditions imply that $S$ is a group.
\noproof\end{theorem}

This result extends the result of \cite{Runde3}, where it
is shown that $M(G)$, the \emph{measure algebra} of a locally
compact group $G$, is Connes-amenable if and only if $G$
is amenable.  This follows as, for discrete groups $G$,
$M(G) = l^1(G)$.

\begin{example}
Let $\omega$ be the weight on $\mathbb Z$ defined by
$\omega(n) = 1+|n|$ for $n\in\mathbb Z$.  By Theorem~\ref{B_AR},
$\mc A = l^1(\mathbb Z,\omega)$ is Arens regular.  For
$m,n\in\mathbb Z$ and $f = (a_k)_{k\in\mathbb Z}\in
l^\infty(\mathbb Z)$, we have that
\[ \ip{ \delta_m \cdot f }{ \delta_n } =
\ip{f}{\delta_{n+m}\Omega(n,m)}
= f_{n+m} \frac{1+|n+m|}{(1+|n|)(1+|m|)}. \]
Suppose that $M \aone \kappa_{\mc A}(\delta_m) = \kappa_{\mc A}(a)$
for some $m\in\mathbb Z$, $M \in l^\infty(\mathbb Z)'$ and $a\in\mc A$.
Then $\ip{M}{\delta_m\cdot f} = \ip{f}{a}$ for each
$f\in l^\infty(\mathbb Z)$, so by letting $f = \kappa_{c_0(\mathbb Z)}(e_k)
\in c_0(\mathbb Z)$, we see that $a = \sum_{k\in\mathbb Z} a_k \delta_k$,
where $a_k = \ip{M}{\delta_m\cdot\kappa_{c_0(\mathbb Z)}(e_k)}$.
However, $\delta_m\cdot\kappa_{c_0(\mathbb Z)}(e_k) \in
\kappa_{c_0(\mathbb Z)}(c_0(\mathbb Z))$ for each $k\in\mathbb Z$,
so if $M \in c_0(\mathbb Z)^\circ$, then $a=0$.

Consequently, if $M \aone \kappa_{\mc A}(\delta_m) \in \kappa_{\mc A}(\mc A)$
for each $m\in\mathbb Z$ and $M \in l^\infty(\mathbb Z)'$, then
$\delta_m \cdot f \in \kappa_{c_0(\mathbb Z)}(c_0(\mathbb Z))$ for
each $m\in\mathbb Z$ and $f \in l^\infty(\mathbb Z)$.  However, if
$\mathbf{1} \in l^\infty(\mathbb Z)$ is the constant $1$ sequence,
then
\[ \lim_n \ip{\delta_m\cdot\mathbf{1}}{\delta_n}
= \lim_n \frac{1+|n+m|}{(1+|n|)(1+|m|)}
= \frac{1}{1+|m|}, \]
so that $\delta_m\cdot\mathbf{1} \not\in
\kappa_{c_0(\mathbb Z)}(c_0(\mathbb Z))$.

We hence conclude that $\mc A$ is not an ideal in $\mc A''$,
and so we cannot apply Theorem~\ref{ca_facts} in this case.
\noproof\end{example}

Unfortunately, it is not possible for $l^1(S,\omega)$ to
be both amenable and Arens regular.

\begin{theorem}\label{Gron_Thm}
Let $G$ be discrete group, and let $\omega$ be a weight on $G$.
Then $l^1(G,\omega)$ is amenable if and only if $G$ is
an amenable group, and $\sup\{ \omega(g) \omega(g^{-1})
: g\in G \} < \infty$.
\end{theorem}
\begin{proof}
This is \cite[Theorem~3.2]{Gron1}.
\end{proof}

\begin{proposition}
Let $S$ be a discrete, unital semigroup, and let $\omega$ be a weight
on $S$ such that $\mc A = l^1(S,\omega)$ is Arens regular.
Let $K>0$ and $B\subseteq S$ be such that for each $g\in B$, $g$
has a right inverse $g^{-1}$ (which need not be unique),
and $\omega(g)\omega(g^{-1})\leq K$.  Then $B$ is finite.
\end{proposition}
\begin{proof}
For $g\in B$ and $h\in S$, we have
\[ \omega(g)\omega(h) = \omega(g) \omega(hgg^{-1}) \leq
\omega(g) \omega(hg) \omega(g^{-1}) \leq K \omega(hg), \]
so that $\Omega(h,g) \geq K^{-1}$.  Suppose now that $B$ is infinite.
Then we can easily construct sequences which violate condition (2)
of Theorem~\ref{B_AR}, showing that $\mc A$ is not Arens regular.
This contradiction shows that $B$ must be finite. 
\end{proof}

\subsection{Injectivity of the predual module}

Let $S$ be a unital, weakly cancellative semigroup, let $\omega$ be a
weight on $S$, and let $\mc A=l^1(S,\omega)$, $\mc A_* = c_0(S)$.
Then $\mc B(\mc A,\mc A_*) = \mc B(l^1,c_0) = l^\infty(c_0) \subseteq
l^\infty(S\times S)$, where we identify $T:\mc A\rightarrow\mc A_*$ with
the bounded family $(\ip{\delta_s}{T(\delta_t)})_{(s,t)\in S\times S}$.
Let $\phi : \mc B(\mc A,\mc A_*) \rightarrow \mc A_*$, so that
$\phi$ is represented by a bounded family $(M_s)_{s\in S} \subseteq
\mc B(\mc A,\mc A_*)'$ using the relation
\[ \ip{\delta_s}{\phi(T)} = \ip{M_s}{T}
\qquad (s\in S, T\in\mc B(\mc A,\mc A_*)). \]
Suppose further that $\phi$ is a left $\mc A$-module homomorphism.
Then
\begin{equation}
\ip{\delta_s}{\phi(T)} = \ip{\delta_{u_s}}{\phi(\delta_s\cdot T)}
= \ip{M_{u_S}}{\delta_s\cdot T} = \ip{M_s}{T}
\qquad (s\in S, T\in\mc B(\mc A,\mc A_*)), \label{eq:three}
\end{equation}
so that $M_s = M_{u_S} \cdot \delta_s$ for each $s\in S$.
We see also that $\phi$ maps into $c_0(S)$ (and not just $l^\infty(S)$)
if and only if
\[ \lim_{s\rightarrow\infty} \ip{M_{u_S}}{\delta_s\cdot T}
= 0 \qquad (T\in\mc B(\mc A,\mc A_*)). \]

Conversely, if condition (\ref{eq:three}) holds, then for $s,t\in S$ and
$T\in\mc B(\mc A,\mc A_*)$, we have that
\begin{align*}
\ip{\delta_s}{\phi(\delta_t\cdot T)} &= \ip{M_s}{\delta_t\cdot T}
= \ip{M_{u_S}}{\delta_s\cdot\delta_t\cdot T}
= \Omega(s,t) \ip{M_{st}}{T} \\
&= \Omega(s,t)\ip{\delta_{st}}{\phi(T)}
= \ip{\delta_s}{\delta_t\cdot \phi(T)}.
\end{align*}
Hence $\phi$ is a left $\mc A$-module homomorphism.

Notice that $c_0(S\times S) \subseteq \mc B(\mc A,\mc A_*)$,
so that $c_0(S\times S)^\circ \subseteq \mc B(\mc A,\mc A_*)'$.

\begin{definition}
Let $G$ be a group and $\omega$ be a weight on $G$ such that for
each $\epsilon>0$, the set $\{ g\in G : \omega(g) \omega(g^{-1})
< \epsilon^{-1} \}$ is finite.  Then we say that the weight
$\omega$ is \emph{strongly non-amenable}.
\end{definition}

\begin{proposition}
Let $G$ be a group, and let $\omega$ be a weight on $G$ such that
$\omega$ is not strongly non-amenable, and let
$\phi:\mc B(\mc A,c_0(G))\rightarrow c_0(G)$ be a left $\mc A$-module
homomorphism.  If $\phi$ is represented by $(M_g)_{g\in G}$ as above,
then $M_{u_G} \in c_0(S\times S)^\circ$.
\end{proposition}
\begin{proof}
We adapt the methods of \cite{DP} to the weighted, discrete case.
As $\omega$ is not strongly non-amenable, there exists some $K>0$
such that the set $X_K = \{ g\in G : \omega(g)\omega(g^{-1})\leq K \}$
is infinite.
Let $M=M_{u_G}$, and suppose that $M\not\in c_0(G\times G)^\circ$,
so that for some $g,h\in G$, we have that $\delta:=\ip{M}{e_{(g,h)}}\not=0$.
We shall henceforth treat $e_{(g,h)}$ as a member of $\mc B(\mc A,c_0(G))$,
noting that for $k\in G$,
\[ \ip{\delta_s}{(\delta_k\cdot e_{(g,h)})(\delta_t)}
= \begin{cases} \Omega(t,k) &: s=g, t=hk^{-1}, \\
0 &: \text{otherwise.} \end{cases} \]
We claim that we can find a sequence $(g_n)_{n\in\mathbb N}$
of distinct elements in $G$ such that
\begin{gather*}
|\ip{M\cdot\delta_{g_m^{-1}g_n}}{e_{(g,h)}}| \leq
   K^{-1} 2^{-2-|m-n|}   \qquad (n\not=m), \\
\omega(g_n) \omega(g_n^{-1}) \leq K \qquad (n\in\mathbb N).
\end{gather*}
We can do this as $\phi$ must map into $c_0(G)$, so that for
any $T:\mc A\rightarrow c_0(G)$, we have $\lim_{g\rightarrow\infty}
\ip{M\cdot\delta_g}{T}=0$.  Explicitly, let $g_1\in X_K$ be arbitrary,
and suppose that we have found $g_1,\ldots,g_k$.  Then notice that
the sets
\begin{gather*}
\{ s\in G : |\ip{M\cdot\delta_{s^{-1}g_n}}{e_{(g,h)}}| > K^{-1}
   2^{-2-|k+1-n|} : 1\leq n\leq k \}, \\
\{ s\in G  : |\ip{M\cdot\delta_{g_m^{-1}s}}{e_{(g,h)}}| > K^{-1}
   2^{-2-|k+1-m|} : 1\leq m\leq k \}
\end{gather*}
are finite, so as $X_K$ is infinite, we can certainly find
some $x_{k+1}$.

Then, for $x=(x_n)\in l^\infty(\mathbb N)$, define $T_x:
\mc A\rightarrow c_0(G)$ by setting
$\ip{\delta_g}{T_x(\delta_{hg_n^{-1}})} = x_n \Omega(hg_n^{-1},g_n)$
for $n\geq 1$, and $\ip{\delta_s}{T_x(\delta_t)}=0$ otherwise.
Then clearly $T_x$ does map into $c_0(G)$, and $\|T_x\| \leq 
\|x\|$.  Notice that for $s,t\in G$, we have
\begin{align*}
\ip{\delta_s}{T_x(\delta_t)}
&= \begin{cases} x_n \Omega(t,g_n) &:
   s=g, t=hg_n^{-1}, \\ 0 &: \text{otherwise,} \end{cases} \\
&= \sum_n x_n \ip{\delta_s}{(\delta_{g_n}\cdot e_{(g,h)})(\delta_t)}.
\end{align*}
Define $Q:l^\infty(\mathbb N)\rightarrow c_0(\mathbb N)$ by
\[ \ip{\delta_n}{Q(x)} = \ip{M}{\delta_{g_n^{-1}}\cdot T_x}
\qquad (n\in\mathbb N), \]
so that $Q$ is bounded and linear.

Let $n_0\geq 1$ and let $x = e_{n_0} \in c_0(\mathbb N)
\subseteq l^\infty(\mathbb N)$.  Then, $T_x = \delta_{g_{n_0}}
\cdot e_{(g,h)}$, so that
\begin{align*}
\ip{\delta_n}{Q(x)} &= \ip{M}{\delta_{g_n^{-1}}\cdot T_x}
= \ip{M}{\delta_{g_n^{-1}}\cdot(\delta_{g_{n_0}}\cdot e_{(g,h)})} \\
&= \begin{cases} \delta\, \Omega(g_{n_0}^{-1},g_{n_0}) &: n=n_0, \\
\Omega(g_n^{-1},g_{n_0})
\ip{M\cdot e_{g_n^{-1}g_{n_0}}}{e_{(g,h)}} &: n\not=n_0. \end{cases}
\end{align*}
Define $Q_1 \in\mc B(c_0(\mathbb N))$ by
\[ Q_1(x) = \big( \Omega(g_n^{-1},g_n) x_n \big)_{n\in\mathbb N}
\qquad ( x=(x_n) \in c_0(\mathbb N) ). \]
Then, as each $g_n\in X_K$, $Q_1$ is an invertible operator.
Let $Q_2$ be the restriction of $Q$ to $c_0(\mathbb N)$,
so that $Q_2\in\mc B(c_0(\mathbb N))$ and $Q_2 = \delta Q_1 +
\delta Q_3Q_1$ for some $Q_3\in\mc B(c_0(\mathbb N))$.
Thus $Q_3 = \delta^{-1}Q_2Q_1^{-1} - I_{c_0(\mathbb N)}$, so that
for $x\in c_0(\mathbb N)$, we have that
\begin{align*}
\|Q_3(x)\| &= \sup_n |\ip{\delta_n}{\delta^{-1}Q_2Q_1^{-1}(x) - x}|
= \sup_n \Big| \sum_m x_m \ip{\delta_n}{\delta^{-1}Q_2Q_1^{-1}(e_m) - e_m} \Big| \\
&= \sup_n \Big| \sum_{m\not=n} x_m \Omega(g_m^{-1},g_m)^{-1}
   \Omega(g_n^{-1},g_m) \ip{M\cdot\delta_{g_n^{-1}g_m}}{e_{(g,h)}} \Big| \\
&\leq K^{-1} \sup_n \sum_{m\not=n} |x_m| 2^{-2-|m-n|} \omega(g_m)\omega(g_m^{-1})
\leq \|x\| / 2.
\end{align*}
Consequently $Q_3-I_{c_0(\mathbb N)}$ is invertible, so that
$Q_2Q_1^{-1}$ is invertible, showing that $Q_2$ is invertible.
However, this implies that $Q_2^{-1}Q:l^\infty(\mathbb N)\rightarrow
c_0(\mathbb N)$ is a projection, which is a well-known contradiction,
completing the proof.
\end{proof}

\begin{theorem}\label{amen_not_inj}
Let $G$ be a countable group, let $\omega$ be a weight which is
not strongly non-amenable, and let $\mc A=l^1(G,\omega)$.
Then $c_0(G)$ is not left-injective.
\end{theorem}
\begin{proof}
Suppose, towards a contradiction, that
$c_0(G)$ is left-injective, so that there exists $M=M_{u_G}
\in \mc B(\mc A,\mc A_*)'$ as above, with the additional condition that
\begin{align*}
\delta_{g,h} &= \ip{\delta_g}{\phi\Delta'_{\mc A}(e_h)}
= \ip{M}{\delta_g \cdot \Delta'_{\mc A}(e_h)}
= \Omega(hg^{-1},g) \ip{M}{\Delta'_{\mc A}(e_{hg^{-1}})} \\
&= \Omega(hg^{-1},g) \ip{M}{\big( \delta_{st,hg^{-1}}
   \Omega(s,t) \big)_{(s,t)\in G\times G}}
\qquad (g,h\in G). \end{align*}
This clearly reduces to
\[ \delta_{g,u_G} = \ip{M}{\big( \delta_{st,g}
   \Omega(s,t) \big)_{(s,t)\in G\times G}} \qquad (g\in G). \]

As $G$ is countable, we can enumerate $G$ as $G = \{ g_n : n\in\mathbb N\}$.
Then, for $g_n\in G$, let $X_{g_n} = \{ g_1, \ldots, g_n \}
\subseteq G$.  Define $Q:l^\infty(G)\rightarrow\mc B(\mc A,c_0(G))$ by
\[ \ip{\delta_s}{Q(x)(\delta_t)} =
\Omega(s,t) \sum_{g\in X_t} x_g \delta_{st,g}
\qquad (s,t\in G, x\in l^\infty(G)). \]
Then, for each $t\in G$, as $X_t$ is finite, we see that
$Q(x)(\delta_t)\in c_0(G)$, so $Q$ is well-defined.
Clearly $Q$ is linear, and we see that for $x\in l^\infty(G)$,
\[ \|Q(x)\| = \sup_{s,t\in G} \Omega(s,t) \Big|\sum_{g\in X_t} x_g \delta_{st,g}\Big|
\leq \sup_{s,t\in G} \sum_{\{g\in X_t : g=st\}} |x_g|
= \|x\|, \]
so that $Q$ is norm-decreasing.  Then, for $h\in G$, we have that
\[ \ip{\delta_s}{Q(e_h)(\delta_t)}
= \Omega(s,t) \sum_{g\in X_t} \delta_{g,h} \delta_{st,g}
= \begin{cases} \ip{\delta_s}{\Delta_{\mc A}'(e_h)(\delta_t)}
&: h\in X_t, \\ 0 &: h\not\in X_t. \end{cases} \]
Let $h=g_{n_0}$, so that $\{ t\in G : h\not\in X_t \} =
\{ g_n\in G : h\not\in X_{g_n} \} = \{ g_1, g_2, \ldots, g_{n_0-1}\}$.
We hence see that $Q(e_{g_0}) - \Delta'_{\mc A}(e_{g_0}) \in
c_0(G\times G)$.  By the preceding proposition, we hence have that
$I_{c_0(G)} = \phi\circ\Delta_{\mc A}' = \phi\circ (Q|_{c_0(G)})$.
However, this implies that $\phi\circ Q:l^\infty(G)\rightarrow
c_0(G)$ is a projection onto $c_0(G)$, giving us the required contradiction.
\end{proof}

\begin{theorem}\label{sg_not_inj}
Let $S$ be a discrete, weakly cancellative semigroup, let $\omega$
be a weight on $S$, and let $\mc A=l^1(S,\omega)$.  When $S$ is
unital, or $S$ is cancellative, $c_0(G)$
is not a bi-injective $\mc A$-bimodule.
\end{theorem}
\begin{proof}
Suppose, towards a contradiction, that $c_0(G)$ is bi-injective.
Then $\mc A$ is Connes-amenable, so that Theorem~\ref{main_thm}
implies that $\mc A$ is amenable, and that $S=G$ is a group.
By Theorem~\ref{Gron_Thm}, we know that $\omega$ is not strongly
non-amenable.  Suppose
that $G$ is countable, so that the above theorem shows that
$c_0(G)$ is not left-injective, and that hence $c_0(G)$ is certainly
not bi-injective, a contradiction.

Suppose that $G$ is not countable.  Then let $H$ be some countably
infinite subgroup of $G$.  Let $K = \sup\{ \omega(g)\omega(g^{-1})
: g\in G \}<\infty$, and let $g,h\in G$.  Then
\[ \Omega(g,h) = \frac{\omega(gh)}{\omega(g)\omega(h)}
= \frac{\omega(gh)}{\omega(g)\omega(g^{-1}gh)}
\geq \frac{\omega(gh)}{\omega(g)\omega(g^{-1})\omega(gh)}
= \frac{1}{\omega(g)\omega(g^{-1})} \geq K^{-1}, \]
so that $\Omega$ is bounded below on $G\times G$, and hence on
$H\times H$.

Then we can find $X\subseteq G$ such that $G = \bigcup_{x\in X} Hx$
and $Hx\cap Hy = \emptyset$ for distinct $x,y\in X$.  Notice that
if $g\in Hx$ then $g^{-1}\in x^{-1}H$, so that $G = \bigcup_{x\in X}
x^{-1}H$ as well.

By the proof of Theorem~\ref{amen_not_inj}, we see that $c_0(H)$ is not
a left-injective $l^1(H,\omega)$-module.  Suppose, towards a
contradiction, that we do have some left $\mc A$-module homomorphism
$\phi:\mc B(l^1(G,\omega),c_0(G)) \rightarrow c_0(G)$ with $\phi
\Delta_{\mc A}' = I_{\mc A'}$.  Notice that
certainly $\mc B(l^1(G,\omega),c_0(G))$ and $c_0(G)$ are Banach left
$l^1(H,\omega)$-modules, by restricting the action from $l^1(G,\omega)$.

Define a map $\psi : \mc B(l^1(H,\omega),c_0(H)) \rightarrow
\mc B(l^1(G,\omega),c_0(G))$ by, for $g,k\in G$,
\[ \ip{\delta_g}{\psi(T)(\delta_k)} =
\begin{cases} \frac{\omega(s) \omega(t)}{\omega(tx) \omega(k)} \ip{\delta_t}{T(\delta_s)}
&: g=tx, k=x^{-1}s \text{ for some } x\in X, s,t\in H, \\
0 &: \text{otherwise.} \end{cases} \]
Certainly $\psi$ is linear, while
\[ \|\psi(T)\|
\leq \|T\| \sup_{s,t\in H, x\in X}
   \frac{\omega(s)\omega(t)}{\omega(tx)\omega(x^{-1}s)}
\leq \|T\| \sup_{s,t\in H, x\in X}
   \frac{\omega(s)\omega(t)}{\omega(txx^{-1}s)}
= \|T\| \sup_{s,t\in H} \Omega(t,s)^{-1}, \]
so that $\psi$ is bounded.  For $h,s,t\in H$, and $x\in X$, we have
\begin{align*}
&\ip{\delta_{tx}}{(\delta_h\cdot\psi(T))(\delta_{x^{-1}s})}
= \Omega(x^{-1}s,h) \ip{\delta_{tx}}{\psi(T)(\delta_{x^{-1}sh})} \\
&= \frac{\Omega(x^{-1}s,h) \omega(sh) \omega(t)}{\omega(tx)
   \omega(x^{-1}sh)} \ip{\delta_t}{T(\delta_s)}
= \frac{\omega(sh) \omega(t)}{\omega(x^{-1}s) \omega(h) \omega(tx)}
   \ip{\delta_t}{T(\delta_s)} \\
&= \omega(s) \omega(x^{-1}s)^{-1} \omega(t) \omega(tx)^{-1}
   \ip{\delta_t}{(\delta_h\cdot T)(\delta_s)}
= \ip{\delta_{tx}}{\psi(\delta_h\cdot T)(\delta_{x^{-1}s})}.
\end{align*}
Thus $\psi$ is a left $l^1(H,\omega)$-module homomorphism.
For $h,s,t\in H$ and $x\in X$, we then have that 
\begin{align*}
\ip{\delta_{tx}}{\psi(\Delta'_{l^1(H,\omega)}(e_h))(\delta_{x^{-1}s})}
&= \frac{\omega(t) \omega(s)}{\omega(x^{-1}s) \omega(tx)}
   \ip{\delta_t}{\delta_s\cdot e_h}
= \Omega(tx,x^{-1}s) \delta_{ts,h} \\
&= \ip{\delta_{tx}}{\delta_{x^{-1}s} \cdot e_h}
= \ip{\delta_{tx}}{\Delta'_{\mc A}(e_h)(\delta_{x^{-1}s})}.
\end{align*}
If $g,k\in G$ are such that $gk\not\in H$ then $g=tx$ and
$k=y^{-1}s$ for some $s,t\in H$ and distinct $x,y\in X$.  Then,
for $h\in H$, we have that $gk\not=h$, so that
\[ \ip{\delta_g}{\Delta'_{\mc A}(e_h)(\delta_k)}
= \Omega(g,k) \delta_{gk,h}
= 0 = \ip{\delta_g}{\psi(\Delta'_{l^1(H,\omega)}(e_h))(\delta_k)}. \]
Hence $\psi\circ \Delta'_{l^1(H,\omega)}$ is equal to
$\Delta'_{\mc A}$ restricted to $l^1(H,\omega)$.

Let $P : c_0(G) \rightarrow c_0(H)$ be the natural projection, which
is obviously an $l^1(H,\omega)$-module homomorphism.  Then
$Q = P \circ \phi \circ \psi : \mc B( l^1(H,\omega), c_0(H) )
\rightarrow c_0(H)$ is a bounded left $l^1(H,\omega)$-module
homomorphism, and $Q \circ \Delta'_{l^1(H,\omega)} = I_{c_o(H)}$.
This contradiction completes the proof.
\end{proof}

We note that just because $\Omega$ is bounded below does not
imply that $\omega$ is bounded, so that $l^1(G,\omega)$ is not
necessarily isomorphic to $l^1(G)$, and hence we cannot simply
apply the results of \cite{DP}.

We have not been able to establish if $c_0(S)$ can every be a
left-injective $l^1(S,\omega)$-module for some semigroup $S$ and
weight $\omega$.

\section{Open questions}

We state a few open questions of interest:

\begin{enumerate}
\item Let $\mc A$ be an Arens regular Banach algebra such that
   $\mc A''$ is Connes-amenable.  Need $\mc A$ be amenable?
\item This is true for C$^*$-algebras.  Can we find a ``simple'' proof?
\item Let $\mc A$ be a dual Banach algebra with predual $\mc A_*$,
   and suppose that $\mc A_*$ is bi-injective.  If $\mc A$ necessarily
   a von Neumann algebra or the bidual of an Arens regular Banach algebra
   $\mc B$ such that $\mc B$ is an ideal in $\mc A$?
\item Let $S$ be a (weakly cancellative) semigroup, and let $\omega$
   be a weight on $S$.  Classify (up to isomorphism) the preduals of
   $l^1(S,\omega)$, and calculate which preduals yield a Connes-amenable
   Banach algebra.
\item This question was asked by Niels Gr{\o}nb\ae k.  In most of our
   examples, it is obvious that when $\mc A$ is a Connes-amenable
   dual Banach algebra, there is $\mc B\subseteq\mc A$ which is
   weak$^*$-dense and amenable.  Is this always true?
\end{enumerate}

\vfill
\begin{list}{}{}
\item[\emph{Address:}] St. John's College,\\
Oxford\\
OX1 3JP\\
United Kingdom
\item[\emph{E-mail:}] \texttt{matt.daws@cantab.net}
\end{list}

\end{document}